\documentstyle[11pt,leqno,graphicx,amscd,fancyhdr,amssymb,latexsym,amsbsy,hyperref,xypic,mathrsfs,pstricks,color,verbatim]{amsart}

\include{micro}
\def\da{\delta}
\def\gz{\gamma}
\def\ggz{\Gamma}
\def\az{\alpha}

\def\bz{\beta}
\def\ra{\rightarrow}
\def\lra{\longrightarrow}
\def\ma{{\mathcal{A}}}
\def\mc{{\mathcal{C}}}

\def\t{\tau}
\def\be{{\textbf{e}}}
\def\Ind{\mbox{Ind}}
\def\bx{{\textbf{x}}}
\def\by{{\textbf{y}}}\def\bg{{\textbf{g}}}
\def\Hom{\mbox{Hom}}
\def\End{\mbox{End}}\def\Ext{\mbox{Ext}}

\newtheorem{thm}{Theorem}[section]
\newtheorem{lem}[thm]{Lemma}
\newtheorem{cor}[thm]{Corollary}

\newtheorem{prop}[thm]{Proposition}
\newtheorem{rem}[thm]{Remark}

\theoremstyle{definition}
\newtheorem{example}[thm]{Example}

\newtheorem{defn}[thm]{Definition}

\begin{document}

\title{A module-theoretic interpretation of Schiffler's expansion formula}
\pagestyle{fancy}
\fancyhead[CO]{A module-theoretic interpretation of Schiffler's expansion formula}
\fancyhead[RE]{}
\fancyhead[RO]{\thepage}
\fancyhead[LO]{}
\fancyhead[LE]{\thepage}
\fancyhead[CE]{ THOMAS BR\"USTLE AND JIE ZHANG }
\lfoot{}
\rfoot{}
\cfoot{}
\renewcommand{\headrulewidth}{0pt}
\renewcommand{\footrulewidth}{0pt}

\author{Thomas Br\"ustle and Jie Zhang}
\address{D\'{e}partement de Math\'{e}matiques, Universit\'{e} de Sherbrooke, Sherbrooke, Canada, J1K 2R1} \email{thomas.brustle@@usherbrooke.ca and tbruestl@@ubishops.ca}
\address{D\'{e}partement de Math\'{e}matiques, Universit\'{e} de Sherbrooke, Sherbrooke, Canada, J1K 2R1} \email{Jie.Zhang@@usherbrooke.ca}

\date{\today}

\thanks{Thomas Br\"ustle is supported by  NSERC.
Jie Zhang is supported by China Scholarship Council.}

\begin{abstract}
We give a module-theoretic interpretation of Schiffler's expansion formula which is defined combinatorially in terms of complete $(\ggz,\gz)-$paths in order to get the expansion of the cluster variables in the cluster algebra of a marked surface $(S,M)$.  Based on the geometric description of the indecomposable objects of the cluster category of the marked surface $(S,M)$, we show the coincidence of
Schiffler-Thomas' expansion formula and the cluster character defined by Palu.
\end{abstract}

\maketitle

\section{Introduction}
We study in this paper the cluster algebra with principal coefficients of a marked surface without punctures:
Consider a compact connected oriented 2-dimensional bordered  Riemann surface $S$ and a finite set of marked points $M$ lying on the boundary $\partial S$ of $S$ with at least one marked point on each boundary component. The condition $M \subset \partial S$ means that we do not allow the marked surface $(S,M)$ to have punctures.

Cluster algebras were first introduced in \cite{FZ02} by Fomin and Zelevinsky. This is a class of commutative rings designed to develop an algebraic framework for the theory of total positivity in semisimple groups and canonical bases in quantum groups. The theory of cluster algebras has connections with various areas of mathematics like quiver representations, Lie theory, combinatorics, Teichm\"uller theory and poisson geometry.

In \cite{FST}, a cluster algebra $\ma_{(S,M)}$ is associated to the marked surface $(S,M)$ with one of its triangulations $\ggz$, and a bijection $\gz\leftrightarrow x^\ggz_\gz$ is established between the internal arcs in $(S,M)$ and the cluster variables in $\ma_{(S,M)}.$  In particular,
the initial seed of  $\ma_{(S,M)}$  corresponds to the triangulation $\Gamma$ of $(S,M)$ and the matrix $B_\ggz$ of $\ggz$, and the mutation of a seed of $\ma_{(S,M)}$ corresponds to the flip of an arc in the triangulation. The cluster algebra $\ma_{(S,M)}$, defined by iterated mutations, thus is independent of the chosen triangulation $\Gamma$ of $(S,M)$.
However, the expression of $x^\ggz_\gz$  was still unclear in \cite{FST} until Schiffler-Thomas \cite{ST} and Schiffler \cite{S} invented the notion of a complete $(\ggz,\gz)-$path which is a concatenation of arcs in $\ggz,$ and completely described the expression of $x^\ggz_\gz$ combinatorially in terms of
$\gz$ and $\ggz$ for the cluster algebra $\ma^0_{(S,M)}$ without coefficients and the cluster algebra $\ma_{(S,M)}$ with coefficients respectively.

A quiver with potential $(Q_{\ggz},W_{\ggz})$ associated to the matrix $B_\ggz$ has been defined for each triangulation $\Gamma$ of $(S,M)$ in \cite{LF,ABCP}. It has been shown in \cite{ABCP} that the corresponding Jacobian algebra $J(Q_\ggz,W_\ggz)$ is a finite-dimensional string algebra. The string modules  over $J(Q_\ggz,W_\ggz)$ are described in \cite{ABCP} by non-contractible curves in the marked surface $(S,M)$ which are not in $\ggz$, and the band modules are described using closed curves in $S$. For each curve $\gz$ in $(S,M)$, we denote by $w(\ggz,\gz)$ and $M(\ggz,\gz)$ the corresponding string and string module over $J(Q_\ggz,W_\ggz)$, respectively.

We give in this paper a module-theoretic interpretation of Schiffler's expansion formula $E^\ggz_?$, interpreting the definition of complete $(\ggz,\gz)-$paths in terms of string modules. We consider the cluster algebra  $\ma_{(S,M)}$ with principal coefficients, and describe Schiffler's expansion formula not only for internal arcs but also for all curves in $(S,M)$ in terms of representations of $J(Q_\ggz,W_\ggz)$. To do this, we denote by $\Ind(\gz)\in{\mathbb{N}}^n$ a vector determined by the minimal injective resolution of $M(\ggz,\gz)$ and by $\mu_{\bf e}(\ggz,\gz)$  the number of collections of substrings of $w(\ggz,\gz)$ such that the direct sum of the corresponding string modules has dimension vector ${\bf e}.$  Then Schiffler's expansion formula $E^\ggz_?$ can be written as follows:

\newtheorem{theorem}{Theorem}
\begin{theorem}Let $\ggz$ be a triangulation of $(S,M)$, and $\gz$ a curve in the marked surface $(S,M)$. Then
$$E^\ggz_\gz=\sum\limits_{\be}\mu_\be(\ggz,\gz)\ X^{\Ind_\ggz(\gz)+B_\ggz\be^T}Y^{\be}.$$
\end{theorem}

The cluster category of the marked surface $\mc_{(S,M)}$ providing a categorification of the cluster algebra
$\ma_{(S,M)}$ has been defined in \cite{A}. In fact, one can define by \cite{A} the cluster
category $\mc_\ggz=\mc_{(Q^{op}_\ggz,W^{op}_{\ggz})}$ associated to $\ggz$ \mbox{since} the Jacobian algebra $J(Q_\ggz,W_\ggz)$ is finite dimensional. Moreover, it has been proven in \cite{KY} that up to triangle equivalence the category $\mc_\ggz$ does not depend on the triangulation $\ggz$ of $(S,M)$ and is just denoted by $\mc_{(S,M)}$.
It is shown in \cite{A} that
$\mc_{(S,M)}$ is a Hom-finite, Krull-Schmidt, 2-Calabi-Yau triangulated category which admits cluster-tilting objects, and the indecomposable objects in $\mc_{(S,M)}$ have been described in \cite{BZ} by curves or closed curves in $(S,M).$

We also show in this paper the coincidence of Schiffler-Thomas' expansion formula and the cluster character defined in \cite{CC,Palu}. In fact,
the cluster-tilting object $T_\ggz$ in $\mc_{(S,M)}$ induced by a triangulation $\ggz$ of $(S,M)$ yields a cluster character
$$X^{\ggz}_?:=X^{T_\ggz}_?: \Ind\mc_{(S,M)}\lra {\mathbb{Q}}(x_1,\ldots,x_n)$$
which sends the rigid indecomposable objects in $\mc_{(S,M)}$ to the cluster variables in the cluster algebra $\ma^0_{(S,M)}$ without coefficients. By using Poettering's formula  for the Euler-Poincar\'e characteristic of a string module, we can get the following theorem:
 \begin{theorem}Let $\gz$ be a curve in $(S,M)$, then
$$X^\ggz_\gz=E^\ggz_\gz(x_1,\ldots,x_n,1,\ldots,1).$$
\end{theorem}
\bigskip

We would like to point out that the equality between the  formulas inherited from complete $(\ggz,\gz)-$paths  and cluster characters has been observed in \cite{DT} for curves in the marked surface $(S,M)$ when $S$ is an annulus.
Moreover, it is shown there that  complete $(\ggz,\gz)-$paths  and cluster characters diverge when considering closed curves in the surface (they give rise to first kind and second kind Chebyshev polynomials respectively).
\bigskip

The paper is organized as follows:
In section 2 we recall the definition of a cluster algebra with coefficients and review the construction of a cluster algebra given by a marked surface. Section 3 summarizes Schiffler's expansion formula for cluster algebras stemming from a surface triangulation.
In section 4 we give module-theoretic interpretations of the combinatorial concepts arising in section 3, and thus obtain a new formulation of Schiffler's expansion formula. We compare these results in section 5 with (Palu's version of) the Caldero-Chapoton map.

\section{Cluster algebras}
\subsection{Definitions}
Let ${\mathbb{P}}=({\mathbb{P}},\oplus,\cdot)$ be a  \emph{semifield} which means that $({\mathbb{P}},\cdot)$ is an abelian multiplicative group, and $\oplus$ is a binary operation on ${\mathbb{P}}$ which is commutative, associative and distributive (with respect to the multiplication $``\cdot"$ ).
Let ${\mathbb{ZP}}$ be the group ring of the multiplicative group $({\mathbb{P}},\cdot)$, and
${\mathcal {F}}$ be the ambient field of rational functions in $n$ independent variables with coefficients in  $\mathbb{QP}.$
\begin{defn}[\cite{FZ07}]A \emph{seed} in $\mathcal {F}$ is a triple $(B,{\bf x,y})$ where
\begin{itemize}
\item[$\bullet$]the \emph{ exchange matrix} $B=(b_{ij})$ is a skew symmetrizable ${n\times n}$ integer matrix,
\item[$\bullet$]the \emph{cluster} ${\bf x}=(x_1,\ldots,x_n)$ is a free generating set of $\mathcal {F}$ over $\mathbb{QP}$ with $n$ algebraically independent elements (called \emph{ cluster variables}),
\item[$\bullet$]the \emph{coefficients} ${\bf y}=(y_1,\ldots,y_n)$ form an $n-$tuple of elements in $\mathbb{P}$.
\end{itemize}
\end{defn}
Denote $[a]_+:=\mbox{max}(x,0)$ for each integer $a$. A new seed can be obtained from a given one
by mutations:
\begin{defn}[\cite{FZ07}]\label{seed mutation}
  The \emph{mutation }of a given seed $(B,{\bf x,y})$ in direction $k\ (1\leq k\leq n)$, denoted by ${\mu_k(B,\bf {x,y})}$, is a new seed $(B',{\bf x',y'})$ in $\mathcal {F}$ defined as follows
\begin{itemize}
\item[$\bullet$]the exchange matrix $B'=(b'_{ij})$ is given by the  matrix mutation defined as follows:
$$b'_{ij}=\left\{
          \begin{array}{ll}
            -b_{ij} & \hbox{\mbox{if} i=k or j=k;} \\
            b_{ij}+[-b_{ik}]_+ \cdot b_{kj}+b_{ik}\cdot[b_{kj}]_+ & \hbox{{otherwise},}
          \end{array}
        \right.
$$
\item[$\bullet$]the cluster ${\bf x'}=({\bf x}\setminus \{x_k\}) \cup \{x'_k\}$ where $x'_k$ is determined by the \emph{exchange relation:}

$$x'_k = \frac{y_k\prod x_i^{[b_{ik}]_+}+\prod x_i^{[-b_{ik}]_+}}{(y_k\oplus 1)x_k},$$

\item[$\bullet$]the coefficients tuple ${\bf y'}=(y'_1,\ldots,y'_n)$ is obtained by mutation of the coefficients tuple defined as follows
$$y'_j=\left\{
         \begin{array}{ll}
           y^{-1}_k & \hbox{\mbox{if}\ j = k;} \\
           y_jy_k^{[b_{kj}]_+}(y_k\oplus 1)^{-b_{kj}} & \hbox{\mbox{if}\ $j\neq k$.}
         \end{array}
       \right.$$
\end{itemize}
\end{defn}
\begin{rem} The mutations are involutions, that is, $(B,{\bf x,y })=\mu_k(\mu_k(B,{\bf x,y })).$ The matrix mutation can also be defined for an $m\times n$ $( m \geq n)$ matrix.
\end{rem}

The \emph{cluster algebra} ${\mathcal {A}=\mathcal {A}}(B,{\bf x,y})$ associated with a given initial seed $(B, \bf {x,y})$ is the $\mathbb{Z}\mathbb{P}-$subalgebra of the field $\mathcal {F}$ generated by all the cluster variables in the clusters obtained by sequences of mutations from the initial seed $(B,{\bf x,y}).$

We usually consider cluster algebras of {\bf\em geometric type} which means
$${\mathbb{P}}=\mbox{Trop}(x_{n+1},\ldots,x_m) (n\leq m)$$
is a  tropical semifield, that is,  $({\mathbb{P}},\cdot)$ is a free abelian group generated by
$x_{n+1},\ldots, x_m$ with addition $\oplus$ defined as follows
$$\prod_{i=n+1}^{m}x_i^{a_i}\oplus\prod_{i=n+1}^{m}x_i^{b_i}=\prod_{i=n+1}^{m}x^{min(a_i,b_i)}_{i\ .}$$

By definition of the tropical semifield $\mathbb{P}$, each element in the coefficients tuple can be written as
$$y_k=\prod_{i=n+1}^m x_{i}^{b_{ik}}$$
where $b_{ik}$ are integers. To deal with the cluster algebras of geometric type, it is convenient to replace the exchange matrix $B$ in Definition \ref{seed mutation} by an $m\times n$ matrix $\tilde{B}=(b_{ij})$ (\emph{extended exchange matrix}) whose upper part is the $n\times n$ exchange matrix $B$ and whose lower part is an $(m-n)\times n$ matrix induced by the coefficient tuple. Then the matrix mutation $\tilde{B}'=(b'_{ij})=\mu_k(\tilde{B})$ yields the mutation of coefficient tuple in Definition \ref{seed mutation}.
Moreover, the exchange relation in Definition \ref{seed mutation} becomes
$$x'_k=x_k^{-1}(\prod_{i=1}^m x_i^{[b_{ik}]_+}+\prod_{i=1}^m x_i^{[-b_{ik}]_+}).$$

A cluster algebra ${\mathcal {A}=\mathcal {A}}(B, {\bf x, y})$ is said to have principal coefficients if it is of gemetric type with the tropical semifield
$${\mathbb{P}}=\mbox{Trop}(x_{n+1},\ldots,x_m)=\mbox{Trop}(y_{1},\ldots,y_n),$$
that is, $m=2n$ and $y_i=x_{n+i}$ for each $1\leq i\leq n$. In the following, we only consider cluster algebras with principal coefficients.
\medskip

\begin{thm}\cite{FZ02,FZ07}Let $X$ be any cluster variable in the cluster algebra $\ma (B, {\bf x, y})$ with principal coefficients, then it has an expansion as
$$X=\frac{f(x_1,x_2,\ldots,x_n;y_1,y_2,\ldots,y_n)}{x_1^{d_1},\ldots,x_n^{d_n}}$$
where $f\in {\mathbb{Z}}[x_1,x_2,\ldots,x_n;y_1,y_2,\ldots,y_n].$

\end{thm}
For each cluster variable $X=X(x_1,\ldots,x_n,y_1,\ldots,y_n)$ in a cluster algebra ${\mathcal {A}}(B, {\bf x, y})$ with principal coefficients, the \emph{F-polynomial} $F_X=F_X(y_1,\ldots,y_n)$
 associated to $X$ is defined to be the rational function obtained from $X$ by evaluating all $x_i$ $(1\leq i\leq n)$ to 1:
$$F_X(y_1,\ldots,y_n)=X(1,\ldots,1,y_1,\ldots,y_n).$$  Let $\tilde{B}=(b_{ij})$ be the $2n\times n$ initial extended exchange matrix of $\mathcal{A}$, then we define $\hat{y_j}$ for each $1\leq j\leq n$ by
$$\hat{y_j}=\prod_{i=1}^m x_i^{b_{ij}}.$$
Proposition 6.1 in \cite{FZ07} implies that each cluster variable $X$ can be written in the form
$$X=X(x_1,\ldots,x_n,y_1,\ldots,y_n)=R(\hat{y}_1,\ldots,\hat{y}_n)\prod_{i=1}^{n}x_i^{g_i},$$
where $g_i$ are integers and $R(u_1,\ldots,u_n)$ is a rational function in ${\mathbb{Q}}(u_1,\ldots,u_n)$. The \emph{g-vector }$g_X$  associated to the cluster variable $X$ is defined as
$$g_X=(g_1,\ldots,g_n).$$

\subsection{Cluster algebras of a marked surface}\label{cluster algebra of the surface}
Let $(S,M)$ be a marked surface without punctures. By a curve in $(S,M)$, we mean the image of a continuous function $\gz : [0,1] \ra S$ with
$\gz(0),\gz(1)\in M$. A simple curve is one where $\gamma$ is injective, except possibly at the endpoints.
We always consider curves up to homotopy, and for any collection of curves we implicitly assume that their mutual intersections are minimal possible in their respective homotopy classes.
The orientation of $S$ induces an orientation on the boundary components (which are circles), and we assume that each such circle is oriented clockwise. We recall from \cite{FST} the definition of a triangulation:

\begin{defn}An \emph{arc} $\da$ in $(S,M)$ is a simple non-contractible curve in $(S,M)$.
The boundary of $S$ is a disjoint union of circles, which are subdivided by the points in $M$ into boundary segments. We call an arc $\delta$ a \emph{boundary arc} if it is homotopic
to such a boundary segment.
Otherwise, $\delta$ is said to be an \emph{internal arc}. A \emph{triangulation} of $(S,M)$ is a maximal
collection $\ggz$ of arcs that do not intersect except at their endpoints. We call a triangle $\triangle$ in $\ggz$ an \emph{internal triangle} if all edges of $\triangle$ are internal arcs.
\end{defn}

Up to homeomorphism, $(S,M)$ is defined by the genus $g$ of $S$, the number $b$ of boundary components, and
the number $c = |M|$ of marked points. Moreover, the following formula gives the number of internal arcs in a triangulation of $(S,M).$

\begin{prop}[\cite{FST}]\label{number-of-arcs}In each triangulation of $(S,M)$, the number
of internal arcs is
$$n = 6g + 3b + c -6.$$
\end{prop}
Recall from \cite{ABCP,LF} that each triangulation $\ggz$ of $(S,M)$ yields
a quiver with potential $(Q_\ggz,W_\ggz)$:
\begin{itemize}
\item[$(1)$] $Q_{\ggz}=(Q_0,Q_1)$ where the set of vertices $Q_0$ is given by the internal arcs of $\ggz$, and the set of arrows $Q_1$ is defined as follows: Whenever there is a triangle $\vartriangle$ in $\ggz$ containing two
internal arcs $a$ and $b$, then there is an arrow $\rho:a \ra b$ in $Q_1$ if $a$ is a predecessor
of $b$ with respect to clockwise orientation at the joint vertex of $a$ and $b$ in $\vartriangle.$
\bigskip

\item[$(2)$]every internal triangle $\vartriangle$ in $\ggz$ gives rise to an oriented cycle $\alpha_\vartriangle \beta_\vartriangle\gz_\vartriangle$ in
$Q$, unique up to cyclic permutation of the factors $\alpha_\vartriangle, \beta_\vartriangle, \gz_\vartriangle$. We define
$$W_\ggz =\displaystyle\sum_\vartriangle\alpha_\vartriangle \beta_\vartriangle\gz_\vartriangle $$
where the sum runs over all internal triangles $\vartriangle$ of $\ggz$.
\end{itemize}

Let $\ggz=\{\t_1,\ldots,\t_n,\t_{n+1},\ldots,\t_m\}$ be a triangulation of $(S,M)$ with $n$ internal arcs $\tau_1, \ldots , \tau_n$.
We define the {\em signed adjacency matrix} $B_\ggz=(b_{ij})_{n\times n}$ of $\ggz$ by setting
$$b_{ij}=\sharp\{\az\in Q_1\mid\az: \t_j\ra \t_i\}-\sharp\{\bz\in Q_1\mid\bz:\t_i\ra \t_j\}.$$
Remark that $B_\ggz$ is always skew-symmetric, and it is the transpose of the matrix of the quiver $Q_\ggz$.

Let ${\mathcal {A}}_\ggz={\mathcal {A}}(B_\ggz, {\bf x}_\ggz, {\bf y}_\ggz)$ be the cluster algebra with principal coefficients associated to $\ggz$, that is ${\mathcal {A}}_\ggz$ is given by the initial seed $(B_\ggz, {\bf x}_\ggz, {\bf y}_\ggz)$ where the exchange matrix is given by the signed adjacency matrix $B_\ggz$, ${\bf x}_\ggz=\{x_{\t_1},\ldots,x_{\t_n}\}$ is the initial cluster indexed by the internal arcs of $\ggz$, the initial coefficients ${\bf y}_\ggz=\{y_{\t_1},\ldots, y_{\t_n}\}$ are given by the vector of generators of a tropical field ${\mathbb{P}}=$Trop$(y_{\t_1},\ldots,y_{\t_n})$. We denote by ${\mathcal{A}}^0_\ggz={\mathcal{A}}(B_\ggz,\bx_\ggz)$ the cluster algebra given by the initial seed $(B_\ggz,\bx,1).$

\begin{thm}[\cite{FST}]\label{FST2} There is a bijection $\gz\longleftrightarrow x^\ggz_\gz$ between the internal arcs in $(S,M)$ and the cluster variables in ${\mathcal{A}}_\ggz$. In particular, ${x}^\ggz_{\t_i}=x_i$ for $1\leq i\leq n$.
\end{thm}

Fomin-Zelevinsky showed in \cite{FZ02} that each cluster variable $x^\ggz_\gz$ is a Laurent polynomial, hence can be uniquely written as $$x^\ggz_\gz=\frac{f(x_1,\ldots,x_n)}{\prod_{i=1}^{n}x_i^{d_i}},$$
where the right hand side is a reduced fraction with $f\in{\mathbb{ZP}}[x_1,\ldots,x_n]$, $d_i\geq 0.$
Fomin-Shapiro-Thurston in \cite{FST} described a recipe to find the denominator vector $(d_1,\ldots,d_n)$ of $x^\ggz_\gz$ in terms of intersection numbers of $\gz$ and the triangulation $\ggz$. Then Schiffler  gave a direct expansion formula of $x^\ggz_\gz$ in \cite{S} by using so-called complete $(\ggz,\gz)-$paths. This
allows to compute $x^\ggz_\gz$ combinatorially in terms of $\gz$ and $\ggz.$

Note that Schiffler and Thomas in \cite {ST} first gave the expansion formula of $x^\ggz_\gz$ for the cluster algebra ${\mathcal {A}}^0_\ggz$.
\section{Schiffler's expansion formula}
We introduce in this section Schiffler's expansion formula for the curves in $(S,M)$. Note that the expansion formula was first defined for internal arcs in $(S,M)$ in order to
get the expression of the corresponding cluster variable (see Theorem \ref{FST2}).  In fact, the definition
of Schiffler's expansion formula works for all curves in $(S,M).$

For two curves $\gz',\gz$ in $(S,M)$ we denote by $I(\gz',\gz)$ the minimal intersection number of two representatives of the homotopic classes of $\gz'$ and $\gz$.

\subsection{Complete $(\ggz,\gz)-$path}
Let $\gz$ be a curve in $(S,M)$ with $d=\sum_{\gz'\in \ggz}I(\gz',\gz)$ which is allowed to have self-intersections.
We fix an orientation for the curve $\gz$ and denote $p_0=\gz(0)$ and $p_{d+1}=\gz(1)$. Let $\t_{i_1}, \t_{i_2},\ldots, \t_{i_d}$ be the internal arcs of $\ggz$ that intersect $\gz$ at $p_1,\ldots, p_d$ in the fixed orientation of $\gz$, see the following figure.

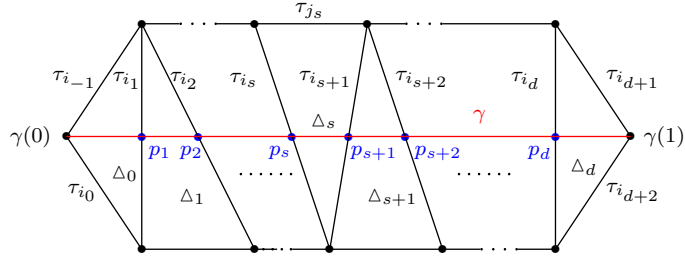
\begin{center}
\begin{pspicture}(0,-1.8)(7.7,1.9)
\put(0,0){\circle*{.1}}\put(7.5,0){\circle*{.1}}
\put(1,1.5){\circle*{.1}}\put(2.5,1.5){\circle*{.1}}\put(4,1.5){\circle*{.1}}
\put(6.5,1.5){\circle*{.1}}
{\blue\put(1,0){\circle*{.12}}\put(1.75,0){\circle*{.12}}
\put(3,0){\circle*{.12}}\put(3.75,0){\circle*{.12}}\put(4.5,0){\circle*{.12}}
\put(6.5,0){\circle*{.12}}
\uput[d](1.25,0.05){$_{p_1}$}\uput[d](1.67,0.05){$_{p_2}$}
\uput[d](2.85,0.05){$_{p_s}$}\uput[d](4.1,0.05){$_{p_{s+1}}$}\uput[d](4.95,0.05){$_{p_{s+2}}$}
\uput[d](6.3,0.05){$_{p_d}$}}

\psline[linewidth=.5pt]{-}(2.8,-1.5)(2.5,-1.5)(1,-1.5)(0,0)(1,1.5)(1.5,1.5)
\psline[linewidth=.5pt]{-}(2,1.5)(2.5,1.5)(4,1.5)(4.5,1.5)
\psline[linewidth=.5pt]{-}(3,-1.5)(3.5,-1.5)(5,-1.5)(5.5,-1.5)
\psline[linewidth=.5pt]{-}(6,-1.5)(6.5,-1.5)
\psline[linewidth=1pt,linestyle=dotted]{-}(5.2,-.5)(5.9,-.5)
\psline[linewidth=1pt,linestyle=dotted]{-}(2.3,-.5)(3,-.5)
\psline[linewidth=.5pt]{-}(6.5,-1.5)(7.5,0)(6.5,1.5)(6.5,-1.5)
\psline[linewidth=.5pt]{-}(1,-1.5)(1,1.5)(2.5,-1.5)\psline[linewidth=.5pt]{-}(2.5,1.5)(3.5,-1.5)(4,1.5)(5,-1.5)
\psline[linewidth=1pt,linestyle=dotted]{-}(1.5,1.5)(2,1.5)\psline[linewidth=1pt,linestyle=dotted]{-}(4.5,1.5)(5,1.5)
\psline[linewidth=.5pt]{-}(6.5,1.5)(5,1.5)
\psline[linewidth=1pt,linestyle=dotted]{-}(2.6,-1.5)(3,-1.5)
\psline[linewidth=1pt,linestyle=dotted]{-}(5.5,-1.5)(6,-1.5)
\uput[l](0,0){$_{\gz(0)}$}\uput[r](7.5,0){$_{\gz(1)}$}
\uput[l](0.6,-.75){$_{\t_{i_{0}}}$}\uput[l](0.6,.75){$_{\t_{i_{-1}}}$}
\uput[u](3.25,1.4){$_{\t_{j_s}}$}
\uput[r](7,.75){$_{\t_{i_{d+1}}}$}\uput[r](7,-.75){$_{\t_{i_{d+2}}}$}
\psline[linecolor=red,linewidth=.5pt](0,0)(7.5,0)
\uput[l](1.2,.75){$_{\t_{i_1}}$}\uput[r](1.2,.75){$_{\t_{i_2}}$}\uput[l](2.75,.75){$_{\t_{i_s}}$}\uput[l](4,.75){$_{\t_{i_{s+1}}}$}
\uput[l](6.5,.75){$_{\t_{i_d}}$}\uput[r](4.2,.75){$_{\t_{i_{s+2}}}$}
\put(.6,-.5){$_{\vartriangle_0}$}\put(1.5,-.8){$_{\vartriangle_1}$}\put(3.25,.2){$_{\vartriangle_s}$}\put(4,-.8){$_{\vartriangle_{s+1}}$}
\put(6.7,-.4){$_{\vartriangle_d}$}{\red \uput[u](5.5,0){$_{\gz}$}}
\put(1,-1.5){\circle*{.1}}\put(2.5,-1.5){\circle*{.1}}\put(3.5,-1.5){\circle*{.1}}\put(5,-1.5){\circle*{.1}}
\put(6.5,-1.5){\circle*{.1}}
\end{pspicture}
\end{center}
Along its way, the curve $\gz$ is passing through (not necessarily distinct) triangles $ \vartriangle_0, \vartriangle_1, \ldots, \vartriangle_d$. For $0\leq k\leq d$, let $[p_kp_{k+1}]$ denote the segment of $\gz$ in $\vartriangle_k$ from $p_k$ to $p_{k+1}$.
If $1\leq l \leq d-1$, then $\vartriangle_l$ is formed by the internal arcs $\t_{i_l}$ and $\t_{i_{l+1}}$, the third arc is denoted by $\t_{j_l}$. In $\vartriangle_0$, we denote the side that lies clockwise
of $\t_{i_1}$ by $\t_{i_{0}}$, the third arc by $\t_{i_{-1}}$. Similarly, in $\vartriangle_d$, denote the side that lies clockwise of $\t_{i_d}$ by $\t_{i_{d+1}}$, the third arc by $\t_{i_{d+2}}$.

\begin{defn}[\cite{ST,S}]\label{complete-path} A {\em complete $(\ggz,\gz)-$path} $\az=\az_1\az_2\cdots \az_{2d+1}$ is a concatenation of arcs in $\ggz$ such that
\begin{itemize}
  \item[$(C1)$] $\az_{2k}=\t_{i_k}$ for each $1\leq k\leq d,$
  \item[$(C2)$] for each $0\leq k\leq d,$ $[p_kp_{k+1}]$ is homotopic to the segment of $\az_{2k}\az_{2k+1}\az_{2k+2}$ in $\az$ starting from $p_k$ to $p_{k+1}.$

\end{itemize}
\end{defn}
Note that a complete $(\ggz,\gz)-$path induces an orientation on each of its arcs $\az_i.$ We denote by $C_\ggz(\gz)$ the set of all complete $(\ggz,\gz)-$paths.
\begin{example}\label{T-r path} Consider a triangulation $\ggz$ of the disc with 8 marked points on the boundary, and let $\gz$ be an internal arc in the disc as follows.
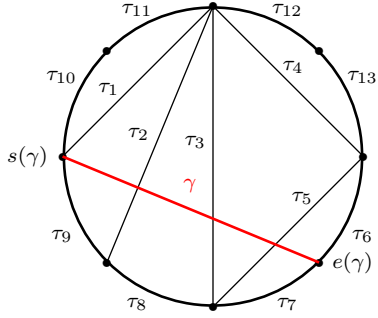
\begin{center}
\begin{pspicture}(-4,-2.2)(4,2.3)
\put(0,-2){\circle*{.1}}\put(-2,0){\circle*{.1}}\put(0,2){\circle*{.1}}\put(2,0){\circle*{.1}}
\put(1.41,1.41){\circle*{.1}}\put(-1.41,1.41){\circle*{.1}}\put(-1.41,-1.41){\circle*{.1}}
\put(1.41,-1.41){\circle*{.1}}
\pscircle[linewidth=1pt](0,0){2}
\psline[linewidth=.5pt](0,2)(2,0)\psline[linewidth=.5pt](0,2)(-2,0)\psline[linewidth=.5pt](0,2)(-1.41,-1.41)
\psline[linewidth=.5pt](0,2)(0,-2)(2,0)
\uput[l](-2,0){$_{s(\gz)}$}\uput[r](1.41,-1.41){$_{e(\gz)}$}
\uput[l](0.1,.2){$_{\t_{3}}$}\uput[r](-1.7,.9){$_{\t_{1}}$}\uput[l](1.4,1.2){$_{\t_{4}}$}
\uput[u](1.2,-.8){$_{\t_{5}}$}\uput[r](-1.3,0.3){$_{\t_{2}}$}
\uput[d](2,-.8){$_{\t_{6}}$}\uput[d](1,-1.7){$_{\t_{7}}$}
\uput[d](-1,-1.7){$_{\t_{8}}$}\uput[d](-2,-.8){$_{\t_{9}}$}
\uput[u](2,.8){$_{\t_{13}}$}\uput[u](1,1.7){$_{\t_{12}}$}
\uput[u](-1,1.7){$_{\t_{11}}$}\uput[u](-2,.8){$_{\t_{10}}$}
\psline[linecolor=red,linewidth=1pt](-2,0)(1.41,-1.41)
{\red\large\uput[u](-.3,-.65){$_\gz$}}
\end{pspicture}
\end{center}
By definition, $C_\ggz(\gz)$ consists of the following 5 elements:
$$\t_1\t_2\t_2\t_3\t_5\t_5\t_7, ~~~ \t_1\t_2\t_8\t_3\t_4\t_5\t_7,$$
$$ \t_1\t_2\t_8\t_3\t_3\t_5\t_6,~~ \t_9\t_2\t_3\t_3\t_3\t_5\t_6,~~~  \t_9\t_2\t_3\t_3\t_4\t_5\t_7.$$
\end{example}

\subsection{Schiffler's expansion formula}\label{schiffler's formula}
Note that the surface $S$ is oriented and the orientation of $S$ induces an orientation on each triangle $\vartriangle$ in $(S,M)$. We assume each triangle $\vartriangle$ has the anticlockwise orientation. Therefore one internal arc $\t_i$ has different orientations in different triangles, see the following picture which is a triangulation of the disk with 8 marked points.
\medskip
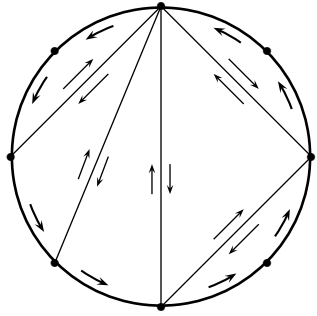
\begin{center}
\begin{pspicture}(-4,-1.7)(4,2)
\put(0,-2){\circle*{.1}}\put(-2,0){\circle*{.1}}\put(0,2){\circle*{.1}}\put(2,0){\circle*{.1}}
\put(1.41,1.41){\circle*{.1}}\put(-1.41,1.41){\circle*{.1}}\put(-1.41,-1.41){\circle*{.1}}
\put(1.41,-1.41){\circle*{.1}}
\pscircle[linewidth=1pt](0,0){2}
\psline[linewidth=.5pt](0,2)(2,0)\psline[linewidth=.5pt](0,2)(-2,0)\psline[linewidth=.5pt](0,2)(-1.41,-1.41)
\psline[linewidth=.5pt](0,2)(0,-2)(2,0)
\psarc[arcsepB=4pt]{->}(0,0){1.85}{55}{70}
\psarc[arcsepB=4pt]{->}(0,0){1.85}{20}{35}
\psarc[arcsepB=4pt]{->}(0,0){1.85}{-70}{-55}
\psarc[arcsepB=4pt]{->}(0,0){1.85}{-35}{-20}
\psarc[arcsepB=4pt]{->}(0,0){1.85}{145}{160}
\psarc[arcsepB=4pt]{->}(0,0){1.85}{110}{125}
\psarc[arcsepB=4pt]{->}(0,0){1.85}{-160}{-145}
\psarc[arcsepB=4pt]{->}(0,0){1.85}{-125}{-110}
\psline[linewidth=.5pt]{->}(.12,-.1)(.12,-.5)\psline[linewidth=.5pt]{<-}(-.12,-.1)(-.12,-.5)
\psline[linewidth=.5pt]{<-}(1.3,.9)(.9,1.3)\psline[linewidth=.5pt]{->}(-1.3,.9)(-.9,1.3)
\psline[linewidth=.5pt]{<-}(1.1,-.7)(.7,-1.1)\psline[linewidth=.5pt]{->}(1.3,-.9)(.9,-1.3)

\psline[linewidth=.5pt]{<-}(-.85,-.4)(-.7,0)\psline[linewidth=.5pt]{->}(-1.1,-.3)(-0.95,0.1)
\psline[linewidth=.5pt]{->}(1.1,.7)(.7,1.1)\psline[linewidth=.5pt]{<-}(-1.1,.7)(-.7,1.1)
\end{pspicture}
\end{center}
\medskip

On the other hand, a complete $(\ggz,\gz)-$path $\az=\az_1\az_2\cdots \az_{2d+1}$ also induces an orientation on each arc $\az_i$ in $\az$. Assume $\az_{2k}=\t_{i_k}$ is a common edge of the two triangles $\vartriangle_{k-1},\vartriangle_k$, then we say $\az_{2k}=\t_{i_k}$ is  \emph{$\gz-$oriented} in $\az$ if the orientation of $\t_{i_k}$ induced by $\az$ is the same as its orientation in  $\vartriangle_k.$

Consider the following segment
$\az_{2k}\az_{2k+1}\az_{2k+2}=\t_{i_k}\t_{j_k}\t_{i_{k+1}}$ of a complete $(\ggz,\gz)-$path $\az=\az_1\az_2\cdots \az_{2d+1}$:
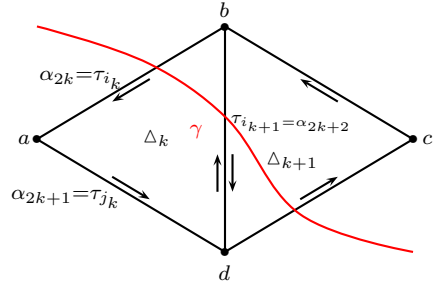
\begin{center}
\begin{pspicture}(-2.5,-1.7)(2.5,1.7)
\put(0,1.5){\circle*{.1}}\put(-2.5,0){\circle*{.1}}\put(2.5,0){\circle*{.1}}\put(0,-1.5){\circle*{.1}}
\psline(0,1.5)(-2.5,0)(0,-1.5)(2.5,0)(0,1.5)(0,-1.5)
\psline{<-}(-.1,-.2)(-.1,-.7)\psline{->}(.1,-.2)(.1,-.7)
\uput[u](0,1.45){$_{b}$}\uput[d](0,-1.5){$_{d}$}\uput[l](-2.4,0){$_{a}$}\uput[r](2.45,0){$_{c}$}
\uput[d](-.9,0.3){$_{\vartriangle_k}$}\uput[d](.9,0){$_{\vartriangle_{k+1}}$}
\psline{<-}(-1.5,.5)(-1,.8)\psline{->}(1.5,.5)(1,.8)
\psline{->}(-1.5,-.5)(-1,-.8)\psline{<-}(1.5,-.5)(1,-.8)
\uput[u](-1.9,.5){$_{\az_{2k}=\t_{i_k}}$}
\uput[d](-2.1,-.5){$_{\az_{2k+1}=\t_{j_k}}$}{\red\uput[l](-.1,.1){$_{\gz}$}}
\pscurve[linecolor=red](-2.5,1.5)(-1,1)(0,.3)(1,-1)(2.5,-1.5)
{\small\uput[r](-.1,.2){$_{\t_{i_{k+1}=\az_{2k+2}}}$}}
\end{pspicture}
\end{center}
Then $\az_{2k}$ is $\gz-$oriented, and $\az_{2k+2}$ is not $\gz-$oriented since the orientation of $\az_{2k+2}$ in $\az$ is from $d$ to $b$ but $\t_{i_{k+1}}=\az_{2k+2}$ is oriented in ${\vartriangle_{k+1}}$ from $b$ towards $d$.

The following proposition given in \cite{S} for internal arcs can also be generalized for curves in $(S,M)$
with the same proof.

\medskip
\begin{prop}[\cite{S}]\label{uniqueness}Let $\gz$ be a curve in $(S,M)$, then there is precisely one complete $(\ggz,\gz)-$path $\az^0_\gz=\az^0_1\az^0_2\cdots \az^0_{2d+1}$ such that none of its even arcs $\az^0_{2k}$ is $\gz-$oriented. Dually, there is precisely one complete $(\ggz,\gz)-$path $\az^1_\gz=\az^1_1\az^1_2\cdots \az^1_{2d+1}$ such that all of its even arcs $\az^1_{2k}$ are $\gz-$oriented.
\end{prop}

\medskip
For each complete $(\ggz,\gz)-$path $\az=\az_1\az_2\cdots \az_{2d+1}$, let
$$x(\az)=\frac{\prod_{k\ odd}x_{\az_k}}{x_{i_1}x_{i_2}\cdots x_{i_d}}\ \mbox{and}\ y(\az)=\prod\limits_{\az_{2k}\ is\ \gz-oriented}y_{i_k}$$
where $x_{\az_k}=1$ if $\az_k$ is a boundary arc.  \emph{Schiffler's expansion formula} $E^\ggz_\gz$  with respect to $\ggz$ is defined as follows:
$$E^\ggz_\gz:=\sum\limits_{\az\in C_\ggz(\gz)}x(\az)y(\az).$$

This formula gives the expression of the cluster variable $x^\ggz_\gz$ in the cluster algebra
${\mathcal{A}}_\ggz:$

\begin{thm}[\cite{S}]Let $\ggz$ be a triangulation of $(S,M)$, and $\gz$ an internal arc in $(S,M)$, then
$$x^\ggz_\gz=E^\ggz_\gz.$$
\end{thm}

\begin{rem}\label{Schiffler-Thomas' formula}If we consider the cluster algebra ${\mathcal{A}}^0_\ggz$, then we get Schiffler-Thomas' expansion formula in \cite{ST}:
$$x^\ggz_\gz=E^\ggz_\gz(x_1,\ldots,x_n,1,\ldots,1)=\sum\limits_{\az\in C_\ggz(\gz)}x(\az)$$
for each internal arc $\gz$.
\end{rem}

\section{A module-theoretic interpretation of Schiffler's expansion formula}
In order to give a module-theoretic interpretation of Schiffler's formula, we first introduce some notation for the string algebra arising from the marked surface $(S,M).$
\subsection{String algebras arising from $(S,M)$}
From \cite{ABCP} we know that  $J(Q_\ggz,W_\ggz)$ is a finite-dimensional string algebra provided $(Q_\ggz,W_\ggz)$ is defined by a triangulation of a marked surface as in the previous section.

In \cite{ABCP},  the strings and bands of $J(Q_\ggz,W_\ggz)$ are related to the curves and closed curves respectively in $(S,M)$ : For each curve $\gz$ in $(S,M)$ with $d=\sum_{\t\in \ggz}I(\t,\gz),$ we fix an orientation of $\gz$, and let $\t_{i_1}, \t_{i_2},\ldots, \t_{i_d}$ be the internal arcs of $\ggz$ that intersect $\gz$ in the fixed orientation of $\gz$.
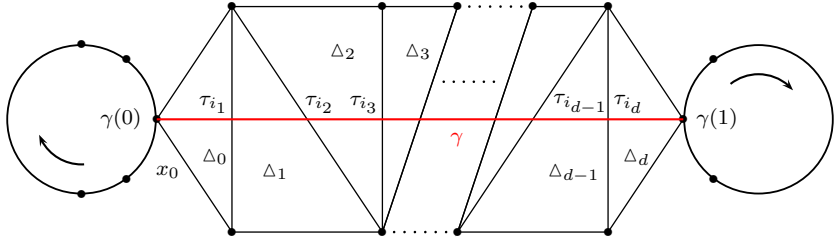
\begin{center}
\begin{pspicture}(-3.4,-2)(12,1.5)
\pscircle[linewidth=.7pt](-2,0){1}
\put(-1.4,.8){\circle*{.1}}\put(-1.4,-.8){\circle*{.1}}\put(-2,1){\circle*{.1}}
\put(-2,-1){\circle*{.1}}\put(-1,0){\circle*{.1}}
\pscircle[linewidth=.7pt](7,0){1}
\put(6.4,.8){\circle*{.1}}\put(6.4,-.8){\circle*{.1}}\put(6,0){\circle*{.1}}
\psarc[arcsepB=4pt]{<-}(-2,0){.6}{-160}{-80}\psarc[arcsepB=4pt]{<-}(7,0){.6}{-320}{-225}
\uput[l](-1,0){$_{\gz(0)}$}\uput[r](6,0){$_{\gz(1)}$}
\put(0,1.5){\circle*{.1}}\put(5,-1.5){\circle*{.1}}\put(5,1.5){\circle*{.1}}
\put(0,-1.5){\circle*{.1}}\put(2,1.5){\circle*{.1}}\put(3,1.5){\circle*{.1}}
\put(4,1.5){\circle*{.1}}\put(2,-1.5){\circle*{.1}}\put(3,-1.5){\circle*{.1}}
\psline[linewidth=.5pt]{-}(-1,0)(0,1.5)(2,1.5)(3,1.5)(2,-1.5)(0,-1.5)(-1,0)
\psline[linewidth=.5pt]{-}(6,0)(5,1.5)(4,1.5)(3,-1.5)(5,-1.5)(6,0)
\psline[linewidth=.5pt]{-}(5,1.5)(3,-1.5)(4,1.5)\psline[linewidth=.5pt]{-}(0,-1.5)(0,1.5)(2,-1.5)(2,1.5)
\psline[linewidth=.5pt]{-}(5,-1.5)(5,1.5)\psline[linewidth=.5pt]{-}(2,-1.5)(3,1.5)
\psline[linestyle=dotted,linewidth=1pt](2,-1.5)(3,-1.5)
\psline[linestyle=dotted,linewidth=1pt](3,1.5)(4,1.5)

\psline[linecolor=red](-1,0)(6,0)
{\red\uput[d](3,0){$_\gz$}}
\psline[linestyle=dotted,linewidth=1pt](3.5,0.5)(2.8,0.5)
\uput[l](0.15,0.2){$_{\t_{i_1}}$}\uput[r](4.9,0.2){$_{\t_{i_d}}$}
\uput[r](.8,0.2){$_{\t_{i_2}}$}\uput[r](4.1,0.2){$_{\t_{i_{d-1}}}$}
\uput[l](2.15,0.2){$_{\t_{i_3}}$}\uput[l](-.5,-0.7){$_{x_0}$}
\put(.4,-0.7){$_{\vartriangle_1}$}\put(-.4,-0.5){$_{\vartriangle_0}$}
\put(1.3,0.9){$_{\vartriangle_2}$}\put(4.2,-0.7){$_{\vartriangle_{d-1}}$}
\put(2.3,0.9){$_{\vartriangle_3}$}\put(5.2,-0.5){$_{\vartriangle_{d}}$}
\end{pspicture}
\end{center}

Here the boundary is indicated by the circles in the above figure.
Along its way, the curve $\gz$ is passing through (not necessarily distinct) triangles $ \vartriangle_0, \vartriangle_1, \ldots, \vartriangle_d$. Thus we obtain a string $w(\ggz,\gz)$ in $J(Q_\ggz,W_\ggz)$:

$$w(\ggz,\gz) : \t_{i_1} \frac{~~~\az_1~}{}\t_{i_2} \frac{~~~\az_2~}{}\cdots \t_{i_{d-2}} \frac{\az_{d-2}}{}\t_{i_{d-1}}\frac{~~~\az_{d-1}~}{}\t_{i_d}.$$
We denote by $M(\ggz,\gz)$ the corresponding string module $M(w(\ggz,\gz))$, we refer to \cite{BR,BZ} for more details about string algebras.

\subsection{Substring Modules}\label{notation}
We keep the notation as in the previous subsection. A {\em subword} $w_{I}(\ggz,\gz)$ of $w(\ggz,\gz)$ indexed by an interval $I=[j,s]=\{j,j+1,\ldots,s\}$ with $1\leq j\leq s\leq d$ is a string over $J(Q_\ggz,W_\ggz)$ given by:
$${w_I(\ggz,\gz)} : \t_{i_j} \frac{~~~\az_j~}{}\t_{i_{j+1}}  \frac{\az_{j+1}}{}\t_{i_{j+2}} \cdots  \frac{~~~\az_{s-1}~}{}\t_{i_{s}}.$$
Note that the string module $M(w_{I}(\ggz,\gz))$ given by the subword $w_{I}(\ggz,\gz)$ of $w(\ggz,\gz)$  is not necessarily a submodule of $M(\ggz,\gz)$. In case $M(w_{I}(\ggz,\gz))$  is  a submodule of $M(\ggz,\gz)$, we call the subword $w_{I}(\ggz,\gz)$  a  \emph{substring} of  $w(\ggz,\gz)$.

More generally, for any subset $I\subset\{1,2,3,\ldots,d\}$, we can uniquely write $I$ as a disjoint union of intervals of maximal length $I=I_1\cup I_2\cup\cdots \cup I_t$ that is
 \begin{itemize}
 \item[(1)]$I_l$ is an interval for each $1\leq l\leq t,$
 \item[(2)]max$\{i|i\in I_l\}+2 \le $min$\{i|i\in I_{l+1}\}$ for each $1\leq l\leq t-1,$
 \item[(3)]$I_j\cap I_k=\emptyset$ if $j\neq k.$
 \end{itemize}
We call $I_1\cup I_2\cup\cdots \cup I_t$ the  \emph{interval decomposition }of $I$.
For each subset $I\subseteq\{1,2,3,\ldots,d\}$ with its interval decomposition $I=I_1\cup I_2\cup\cdots \cup I_t$, consider
the direct sum of string modules over $J(Q_\ggz,W_\ggz)$
$$M_I(\ggz,\gz)=\bigoplus_{l=1}^tM(w_{I_l}(\ggz,\gz)).$$
Set
$$S_\ggz(\gz)=\{I\subseteq\{1,2,3,\ldots,d\}|M_I(\ggz,\gz)\ is\ a\ submodule\ of\ M(\ggz,\gz) \}.$$
Therefore, if $I=I_1\cup I_2\cup\cdots \cup I_t\in S_\ggz(\gz)$ then $w_{I_l}(\ggz,\gz)$ is a substring of $w(\ggz,\gz)$ for each $1\leq l\leq t.$

\begin{example}\label{notation}We consider an annulus with three marked points on the boundary components and a triangulation $\ggz=\{\t_1,\t_2,\t_3,\t_4,\t_5,\t_6\}$ with internal arcs $\t_1,\t_2,\t_3$ as follows:
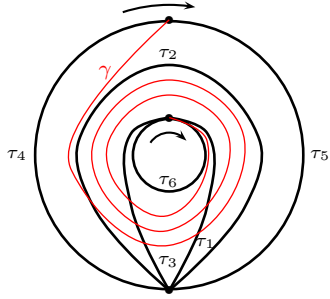
\begin{center}
\begin{pspicture}(-4,-1.8)(4,2)
\put(0,.5){\circle*{.1}}
\pscircle[linewidth=1pt](0,0){.5}
\pscircle[linewidth=1pt](0,0){1.8}
\put(0,1.8){\circle*{.1}}
\put(0,-1.8){\circle*{.1}}
\psarc[arcsepB=4pt]{<-}(0,0){.3}{-320}{-200}\psarc[arcsepB=4pt]{<-}(0,0){2}{80}{110}
\pscurve[linewidth=1pt]{-}(0,-1.8)(.5,-0.7)(.5,.3)(0,.5) \pscurve[linewidth=1pt]{-}(0,-1.8)(-.5,-0.7)(-.5,.3)(0,.5)
\pscurve[linewidth=1pt]{-}(0,-1.8)(1,-0.7)(1.2,.2)(0,1.2)(-1.2,.2)(-1,-0.7)(0,-1.8)
\pscurve[linewidth=.5pt,linecolor=red,showpoints=false]{-}(0,1.8)(-1.3,0.2)(-1.3,-0.2)(0,-1.2)(1,.2)(0,1)(-1,.2)(-.8,-.6)(0,-1)
(.8,.2)(0,.8)(-.8,.2)(-.7,-.4)(0,-.8)(.5,.2)(0,.5)
{\small\uput[u](0,1.1){$_{\t_2}$}\uput[l](-1.7,0){$_{\t_4}$}\uput[r](1.7,0){$_{\t_5}$}
\put(.35,-1.15){$_{\t_1}$}\put(-.15,-1.4){$_{\t_3}$}\uput[u](0,-.6){$_{\t_6}$}}
{\red\uput[l](-.7,1.1){$_\gz$}}
\end{pspicture}
\end{center}
The associated Jacobian algebra $J(Q_\ggz,W_\ggz)$ is cluster-tilted of type $\tilde{A}_2$ 
where the quiver is as follows
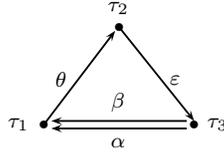
\begin{center}
\begin{pspicture}(-1,-.1)(1,1.3)
\put(0,1.3){\circle*{.1}}\put(-1,0){\circle*{.1}}\put(1,0){\circle*{.1}}
\psline[linewidth=.7pt]{->}(-1,0)(-.05,1.25)\psline[linewidth=.7pt]{<-}(1,0.05)(0,1.3)
\psline[linewidth=.7pt]{->}(.9,0.05)(-.9,0.05)\psline[linewidth=.7pt]{->}(.9,-0.05)(-.9,-0.05)
\uput[u](0,1.3){$_{\t_2}$}\uput[l](-1,0){$_{\t_1}$}\uput[r](1,0){$_{\t_3}$}
\uput[u](0,0){$_{\bz}$}\uput[d](0,0){$_{\az}$}
\uput[l](-.5,.6){$_{\theta}$}\uput[r](.5,.6){$_{\varepsilon}$}
\end{pspicture}
\end{center}
with potential $W_\ggz=\bz\varepsilon\theta.$

For the curve $\gz$, it is easy to see that $d=\sum_{\t\in\ggz}I(\t,\gz)=6$ with $\t_{i_1}=\t_2, \t_{i_2}=\t_{i_4}=\t_{i_6}=\t_3, \t_{i_3}=\t_{i_5}=\t_{1}$ and $w(\ggz,\gz)$ is given by:
$$w(\ggz,\gz): \t_2\overset{\varepsilon}{\ra}\t_3\overset{\az}{\ra}\t_1\overset{\bz}{\leftarrow}\t_3\overset{\az}{\ra}
\t_1\overset{\bz}{\leftarrow}\t_3.$$
Take $I=\{1,3,5\}$, then $\{1\}\cup\{3\}\cup\{5\}$ is the interval decomposition of $I$, and $M_I(\ggz,\gz)=S_2\oplus S_1\oplus S_1$ where $S_i$ is the simple module corresponding to $\t_{i}$. Since $S_2$ is not a submodule of $M(\ggz,\gz)$, $I\not\in S_\ggz(\gz)$.

Take $I=\{1,2,3,5\}=\{1,2,3\}\cup \{5\}$, then $M_I(\ggz,\gz)=M(\ggz,\az\varepsilon)\oplus S_1$, then $I\in S_\ggz(\gz)$ since both $M(\ggz,\az\varepsilon)$ and $S_1$ are submodules of $M(\ggz,\gz).$
Moreover,
$S_\ggz(\gz)$ consists of the following elements
$$\emptyset,~\{3,5\},~ \{1,2,3\},~\{3,4,5\},~ \{1,2,3,5\},~\{2,3,4,5,6\},$$
$$\{3\},\{2,3\},~\{2,3,5\},~ \{2,3,5,6\},~ \{2,3,4,5\},~\{1,2,3,5,6\},$$
$$\{5\},\{5,6\},~\{3,5,6\},~\{3,4,5,6\},~ \{1,2,3,4,5\},~\{1,2,3,4,5,6\}.$$

\end{example}
\subsection{$\gz-$oriented arcs}\label{r-orient}
The definition of a $\gz-$oriented arc in a complete $(\ggz,\gz)-$path is given in combinatorial terms. We give in this section a module-theoretic interpretation for this definition. The following two lemmas play an important role for the interpretation.

\begin{lem}\label{lem1-r-oriented}Let $\az=\az_1\az_2\cdots \az_{2d+1}$ be a complete $(\ggz,\gz)-$path in $(S,M)$. Assume there is an arrow $\az_{2k}=\t_{i_{k}}\ra\t_{i_{k+1}}=\az_{2k+2}$ \emph{(}or $\t_{i_{k-1}}=\az_{2k-2}\leftarrow\t_{i_{k}}=\az_{2k}$\emph{)} in $J(Q_\ggz,W_\ggz)$. If $\az_{2k}$ is $\gz-$oriented in $\az$, then $\az_{{2k+2}}$ \emph{(}or $\az_{2k-2}$\emph{)} is also $\gz-$oriented in $\az$.
\end{lem}
\begin{pf}We only prove the case when there is an arrow $\t_{i_{k}}\ra\t_{i_{k+1}}$ in $J(Q_\ggz,W_\ggz)$.
Assume $\az=\az_1\az_2\cdots \az_{2k-1}\t_{i_k}\az_{2k+1}\t_{i_{k+1}}\az_{2k+3}\cdots\az_{2d+1},$
and consider the following segment $\t_{i_k}\az_{2k+1}\t_{i_{k+1}}$ of $\az$,

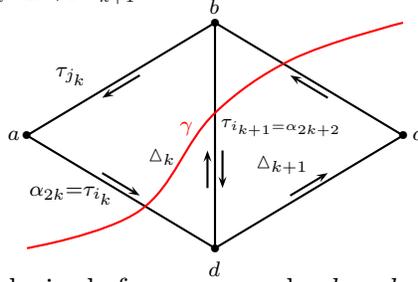
\begin{center}
\begin{pspicture}(-2.5,-1.5)(2.5,1.5)
\put(0,1.5){\circle*{.1}}\put(-2.5,0){\circle*{.1}}\put(2.5,0){\circle*{.1}}\put(0,-1.5){\circle*{.1}}
\psline(0,1.5)(-2.5,0)(0,-1.5)(2.5,0)(0,1.5)(0,-1.5)
\psline{<-}(-.1,-.2)(-.1,-.7)\psline{->}(.1,-.2)(.1,-.7)
\uput[u](0,1.45){$_{b}$}\uput[d](0,-1.5){$_{d}$}\uput[l](-2.4,0){$_{a}$}\uput[r](2.45,0){$_{c}$}
\uput[d](-.7,0){$_{\vartriangle_k}$}\uput[d](.9,-0.1){$_{\vartriangle_{k+1}}$}
\psline{<-}(-1.5,.5)(-1,.8)\psline{->}(1.5,.5)(1,.8)
\psline{->}(-1.5,-.5)(-1,-.8)\psline{<-}(1.5,-.5)(1,-.8)
\uput[d](-1.9,-.5){$_{\az_{2k}=\t_{i_k}}$}{\red\uput[l](-.1,.1){$_{\gz}$}}
\pscurve[linecolor=red](-2.5,-1.5)(-1,-1)(0,.3)(1,1)(2.5,1.5)
{\small\uput[r](-.1,0.1){$_{\t_{i_{k+1}=\az_{2k+2}}}$}}\uput[u](-1.9,.5){$_{\t_{j_k}}$}
\end{pspicture}
\end{center}
 where $\t_{i_k}$ lies anticlockwise before $\t_{i_{k+1}}$, and $a,b,c,d$ are four marked points in $M$.

Since $\t_{i_k}=\az_{2k}$ is $\gz-$oriented in $\az$, $\az_{2k}$ has an orientation from $a$ to $d$ in $\az$. The definition of complete $(\ggz,\gz)-$path implies $\az_{2k+2}=\az_{2k+1}=\t_{i_{k+1}}$ and
 $\az$ induces an orientation for $\az_{2k+2}=\t_{i_{k+1}}$ which is the same as induced by $\vartriangle_{k+1}$. Hence $\t_{i_{k+1}}=\az_{2k+2}$ is  $\gz-$oriented in $\az$.
\end{pf}

The above lemma allows us to define a map
$$\varphi: C_\ggz(\gz)\lra S_\ggz(\gz)$$
by sending each $\az\in C_\ggz(\gz)$ to $I_\az:=\varphi(\az)=\{k|\ \az_{2k}\ \mbox{is}\ \gz-\mbox{oriented} \}$. Note that $M_{I_\az}(\ggz,\gz)$ is a submodule of $M(\ggz,\gz)$ which implies $I_\az\in S_\ggz(\gz).$

\begin{lem}\label{lem2-r-oriented}For each $I\in S_\ggz(\gz)$, there is a unique complete $(\ggz,\gz)-$path
$\az_I=\az_1\az_2\cdots\az_{2d+1}\in C_\ggz(\gz)$ such that only the elements in $\{\az_{2k}|\ k\in I\}$
are $\gz-$oriented.
\end{lem}
\begin{pf}Suppose $I=I_1\cup I_2\cup\cdots \cup I_t$ is the interval decomposition of $I.$ Without loss of generality, we assume for each $1\leq l\leq t$
$$I_l=\{r_l,r_l+1,r_l+2,\ldots,s_l\}\subset\{1,2,3,\ldots,d\}$$
where $1\leq r_1\leq s_1<r_2\leq s_2<\cdots<r_t\leq s_t\leq d.$ Let $\az^0_\gz=\az^0_1\az^0_2\cdots \az^0_{2d+1}$ be the complete $(\ggz,\gz)-$path in Proposition \ref{uniqueness} such that none of its even arcs $\az^0_{2k}$ is $\gz-$oriented and $\az^1_\gz=\az^1_1\az^1_2\cdots \az^1_{2d+1}$ be the one with all its even arcs $\az^1_{2k}$ $\gz-$oriented.

Consider the following expression
$$\az_I=\az^1_0\underline{\az^0_2}\cdots \underline{\az^0_{2r_1-2}}(\az_{{2r_1-1}}\underline{\az^1_{2r_1}}\az^1_{2r_1+1}\cdots\underline{\az^1_{2s_1}}
\az_{{2s_1+1}})$$
$$\underline{\az^0_{2s_1+2}}\cdots\underline{\az^0_{2r_2-2}}(\az_{{2r_2-1}}\underline{\az^1_{2r_2}}\az^1_{2r_2+1}\cdots
\underline{\az^1_{2s_2}}\az_{{2s_2}+1})$$
$$\underline{\az^0_{2s_2+2}}\cdots\az^0_{2r_3-2}\cdots\cdots$$
$$(\az_{{2r_t-1}}\underline{\az^1_{2r_t}}\az^1_{2r_t+1}\cdots
\underline{\az^1_{2s_t}}
\az_{{2s_t+1}})\underline{\az^0_{2s_t+2}}\cdots\az^0_{2d+1}.$$
where the even positions are underlined. We are going to show that
there is a unique choice of $\az_{2r_l-1}$ (or $\az_{2s_l+1}$) for each $1\leq l \leq t$ such that $\az_I$ is a complete
$(\ggz,\gz)-$path.

For each $1\leq l\leq t,$ consider the following segment
$ \az^0_{2r_l-2}\az_{2r_l-1}\az^1_{2r_l}$ of $\az_I,$
\begin{center}
\begin{pspicture}(-2.5,-1.5)(2.5,1.5)
\put(0,1.5){\circle*{.1}}\put(-2.5,0){\circle*{.1}}\put(2.5,0){\circle*{.1}}\put(0,-1.5){\circle*{.1}}
\psline(0,1.5)(-2.5,0)(0,-1.5)(2.5,0)(0,1.5)(0,-1.5)
\psline{<-}(-.1,-.2)(-.1,-.7)\psline{->}(.1,-.2)(.1,-.7)
\uput[u](0,1.45){$_{b}$}\uput[d](0,-1.45){$_{d}$}\uput[l](-2.4,0){$_{a}$}\uput[r](2.45,0){$_{c}$}
\uput[d](-1,.2){$_{\vartriangle_{l-1}}$}\uput[d](.7,-0.1){$_{\vartriangle_{l}}$}
\psline{<-}(-1.5,.5)(-1,.8)\psline{->}(1.5,.5)(1,.8)
\psline{->}(-1.5,-.5)(-1,-.8)\psline{<-}(1.5,-.5)(1,-.8)
\uput[d](-2.1,-.5){$_{\t_{i_{r_l-1}}=\az^0_{2r_l-2}}$}{\red\uput[l](-.1,.1){$_{\gz}$}}
\pscurve[linecolor=red](-2.5,-1.5)(-1,-1)(0,.4)(1,1)(2.5,1.5)
\uput[r](-.1,0.2){$_{\az^1_{2r_l}=\t_{i_{r_l}}}$}
\uput[u](-1.9,.5){$_{\t_{j_{l-1}}}$}
\end{pspicture}
\end{center}
where Lemma \ref{lem1-r-oriented} implies  that $\az^0_{2r_l-2}=\t_{i_{r_l-1}}$ is a clockwise successor of $\az^1_{2r_l}=\t_{i_{r_l}}$.

Since $\az^1_{2r_l}$ is $\gz-$oriented, it has the orientation from $b$ to $d$. And $\az^0_{2r_l-2}$ is oriented from $d$ towards $a$ since $\az^0_{2r_l-2}=\t_{i_{r_l-1}}$ is not $\gz-$oriented in $\az_I.$ Thus $\az_{2r_l-1}=\t_{j_{l-1}}$ for $1\leq l\leq t.$

Similarly, we get $\az_{2s_l+1}=\t_{j_{s_l}}$ for $1\leq l\leq t.$ Hence $\az_I$ is a complete $(\ggz,\gz)-$path such that only the elements in $\{\az_{2k}|\ k\in I\}$ are not $\gz-$oriented in $\az_I.$
The uniqueness of $\az^0$ and $\az^1$ yields the uniqueness of $\az_I.$
\end{pf}
Define a map
$$\psi: S_\ggz(\gz)\lra C_\ggz(\gz)$$
by sending each $I\in S_\ggz(\gz)$ to $\psi(I)=\az_I$ as defined in the above lemma.

It is easy to check that $\psi\circ\varphi (I)=I$ by the definition and $\varphi\circ\psi(\az)=\az$ by the uniqueness
in Lemma \ref{lem2-r-oriented}. Hence $\psi$ and $\varphi$ are mutually inverse, and there is a bijection between $S_\ggz(\gz)$ and $C_\ggz(\gz)$.

\begin{example}Reconsider Example \ref{T-r path}, then
$$\az_\gz^0=\t_1\t_2\t_2\t_3\t_5\t_5\t_7, ~~~\az_\gz^1=\t_9\t_2\t_3\t_3
\t_3\t_5\t_6.$$
The string $w(\ggz,\gz)$ is of the form:
$$\t_2\leftarrow\t_3\ra\t_5,$$
then $S_\ggz(\gz)=\{\emptyset, \{1\}, \{3\}, \{1,3\}, \{1,2,3\}\}.$ Hence,
$$\psi(\emptyset)=\t_9\t_2\t_3\t_3\t_3\t_5\t_6,~~
\psi(\{1,2,3\})=\t_1\t_2\t_2\t_3\t_5\t_5\t_7,$$
$$\psi(\{1\})=\t_1\t_2\t_8\t_3\t_3\t_5\t_6,
~\psi(\{3\})=\t_9\t_2\t_3\t_3\t_4\t_5\t_7,
~\psi(\{1,3\})=\t_1\t_2\t_8\t_3\t_4\t_5\t_7.$$
\end{example}

\subsection{Index of $\gz$}\label{index of r}

Recall from \cite{FZ07} that the Laurent expansion of any cluster variable $x^\ggz_\gz$ in the cluster algebra
${\mathcal {A}}_\ggz$ with initial seed $(B_\ggz, \bx_\ggz,\by_\ggz)$ is homogeneous with respect to the grading given by
$$\deg(x_i)=\deg(x_{\t_i})=\be_{\t_i}=(0,\ldots,0,1,0,\ldots,0)^T$$
with 1 at $i^{th}$ position and
$$\deg(y_i)=-B_\ggz\be_{\t_i}=-B_\ggz\be_{i}$$ where we often write $\be_i$
instead of $\be_{\t_i}$ for convenience. The definition of $B_\ggz$ implies
$$\deg(y_i)=\sum_{\t_i\ra\t_k}\deg(x_k)-\sum_{\t_k\ra\t_i}\deg(x_k)
=\deg(\frac{\prod_{\t_i\ra\t_k}x_k}{\prod_{\t_k\ra\t_i}x_k}).$$

The $\textbf{g}-$vector $\textbf{g}^\ggz_\gz$ associated to the cluster variable $x^\ggz_\gz$ is just the degree of its Laurent expansion with respect to this grading.

Let $\gz$ be a curve in $(S,M)$. If $\gz=\t_i\in\ggz$, define the {\bf index} of $\gz$ with respect to $\ggz$ by $\mbox{Ind}_\ggz(\gz)=\be_i;$
If $\gz\not\in\ggz$ and
  $$ 0\lra M(\ggz,\gz)\lra I_0 \lra I_1 $$
  is the minimal injective resolution of $M(\ggz,\gz)$ in {\bf mod}$J(Q_\ggz,W_\ggz)$, then the {\bf index} ${\Ind_\ggz(\gz)}$ of $\gz$ with respect to $\ggz$ is defined to be the vector given by
  $$\mbox{Ind}_\ggz(\gz)=\sum^n_{i=1}{\bf e_i}~(\mid \Hom(S_i, I_1)\mid-\mid \Hom(S_i,I_0)\mid)$$
  where $S_i$ is the simple module over $J(Q_\ggz,W_\ggz).$

\begin{prop}\label{g-vector}Let $\gz$ be a curve in $(S,M)$, $\az^0_\gz\in C_\ggz(\gz)$ be the complete $(\ggz,\gz)-$path such that none of its even arcs is $\gz-$oriented, then
\begin{itemize}
  \item[$(1)$]$\deg(x(\az_\gz^0))=\emph{Ind}_\ggz(\gz)$,
  \item[$(2)$]$E^\ggz_\gz$ is homogeneous.
\end{itemize}
\end{prop}
\begin{pf}It suffices to consider $\gz\not\in\ggz$, we assume $\gz$ intersects $\ggz$ in its fixed orientation at $\t_{i_1},\ldots,\t_{i_d}$, see the following picture where we write $i_k (\mbox{or}\ j_k )$ instead of $\t_{i_k} (\mbox{or}\ \t_{ j_k } )$ for convenience.
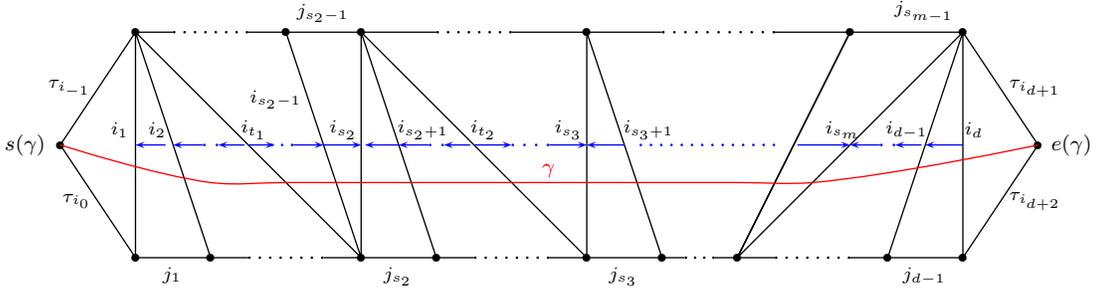
\begin{center}
\begin{pspicture}(0,-1.6)(14,1.5)
\put(0,0){\circle*{.1}}\put(13,0){\circle*{.1}}
\put(1,1.5){\circle*{.1}}\put(3,1.5){\circle*{.1}}
\put(4,1.5){\circle*{.1}}\put(7,1.5){\circle*{.1}}\put(10.5,1.5){\circle*{.1}}
\put(12,1.5){\circle*{.1}}
\put(1,-1.5){\circle*{.1}}\put(4,-1.5){\circle*{.1}}
\put(2,-1.5){\circle*{.1}}
\put(8,-1.5){\circle*{.1}}
\put(7,-1.5){\circle*{.1}}\put(5,-1.5){\circle*{.1}}
\put(9,-1.5){\circle*{.1}}\put(11,-1.5){\circle*{.1}}\put(12,-1.5){\circle*{.1}}

\psline[linewidth=1pt,linestyle=dotted](2.5,1.5)(1.5,1.5)
\psline[linewidth=1pt,linestyle=dotted](3.5,-1.5)(2.5,-1.5)
\psline[linewidth=.5pt]{-}(1,-1.5)(0,0)(1,1.5)(1,-1.5)(2,-1.5)(1,1.5)(4,-1.5)(3.5,-1.5)
\psline[linewidth=.5pt]{-}(4,1.5)(4,-1.5)(3,1.5)(4,1.5)(5,1.5)
\psline[linewidth=.5pt]{-}(2,-1.5)(2.5,-1.5)
\psline[linewidth=.5pt]{-}(1,1.5)(1.5,1.5)\psline[linewidth=.5pt]{-}(2.5,1.5)(3,1.5)
\psline[linewidth=.5pt]{-}(4,-1.5)(5,-1.5)(5.5,-1.5)
\psline[linewidth=.5pt]{-}(5,-1.5)(4,1.5)(7,-1.5)(6.5,-1.5)
\psline[linewidth=1pt,linestyle=dotted](5.5,-1.5)(6.5,-1.5)
\psline[linewidth=1pt,linestyle=dotted](5,1.5)(6,1.5)
\psline[linewidth=.5pt]{-}(6,1.5)(7,1.5)(7,-1.5)(8,-1.5)(7,1.5)(8,1.5)
\psline[linewidth=1pt,linestyle=dotted](8.3,-1.5)(8.7,-1.5)
\psline[linewidth=.5pt]{-}(8,-1.5)(8.3,-1.5)
\psline[linewidth=.5pt]{-}(8.7,-1.5)(9,-1.5)(10.5,1.5)(12,1.5)(9,-1.5)(10.5,1.5)(9.5,1.5)
\psline[linewidth=1pt,linestyle=dotted](8,1.5)(9.5,1.5)
\psline[linewidth=.5pt]{-}(9.5,-1.5)(9,-1.5)
\psline[linewidth=.5pt]{-}(10.5,-1.5)(11,-1.5)(12,1.5)(12,-1.5)(11,-1.5)
\psline[linewidth=.5pt]{-}(12,-1.5)(13,0)(12,1.5)
\psline[linewidth=1pt,linestyle=dotted]{-}(9.5,-1.5)(10.5,-1.5)
\uput[l](0,0){$_{s(\gz)}$}\uput[r](13,0){$_{e(\gz)}$}
{\red \uput[u](6.5,-.6){$_{\gz}$}}
{\small\uput[r](.5,0.2){$_{{i_1}}$}
\uput[l](2.95,0.2){$_{{i_{t_1}}}$}
\uput[r](1,0.2){$_{{i_2}}$}\uput[d](1.5,-1.45){$_{{j_1}}$}
\uput[l](4.15,0.2){$_{{i_{s_2}}}$}\uput[u](3.5,1.45){$_{{j_{s_2-1}}}$}\uput[d](4.5,-1.45){$_{{j_{s_2}}}$}
\uput[l](5.35,0.2){$_{{i_{s_2+1}}}$}\uput[l](3.4,0.6){$_{{i_{s_2-1}}}$}
\uput[l](8.35,0.2){$_{{i_{s_3+1}}}$}
\uput[r](10,0.2){$_{{i_{s_m}}}$}\uput[u](11.5,1.45){$_{{j_{s_{m-1}}}}$}
\uput[l](6,0.2){$_{{i_{t_2}}}$}
\uput[r](6.4,0.2){$_{{i_{s_3}}}$}\uput[d](7.5,-1.45){$_{{j_{s_3}}}$}
\uput[r](11.85,0.2){$_{{i_d}}$}
\uput[r](10.8,0.2){$_{{i_{d-1}}}$}\uput[d](11.5,-1.45){$_{{j_{d-1}}}$}
\uput[l](0.6,-.75){$_{\t_{i_{0}}}$}
\uput[l](0.6,.75){$_{\t_{i_{-1}}}$}
\uput[l](13.5,-.75){$_{\t_{i_{d+2}}}$}
\uput[l](13.5,.75){$_{\t_{i_{d+1}}}$}}

\pscurve[linewidth=.5pt,linecolor=red]{-}(0,0)(2,-.5)(3,-.5)(7,-.5)(9,-.5)(10,-.5)(13,0)
\psline[linewidth=.5pt,linecolor=blue]{->}(1.4,0)(1.0,0)
\psline[linewidth=.5pt,linecolor=blue]{->}(1.9,0)(1.5,0)
\psline[linewidth=.5pt,linecolor=blue]{->}(2.45,0)(2.1,0)
\psline[linewidth=1pt,linestyle=dotted,linecolor=blue]{-}(2.1,0)(1.9,0)
\psline[linewidth=.5pt,linecolor=blue]{<-}(2.85,0)(2.5,0)
\psline[linewidth=.5pt,linecolor=blue]{<-}(3.5,0)(3.1,0)
\psline[linewidth=1pt,linestyle=dotted,linecolor=blue]{-}(3.1,0)(2.8,0)
\psline[linewidth=.5pt,linecolor=blue]{<-}(4,0)(3.5,0)
\psline[linewidth=.5pt,linecolor=blue]{<-}(4.05,0)(4.5,0)
\psline[linewidth=.5pt,linecolor=blue]{->}(4.9,0)(4.5,0)
\psline[linewidth=.5pt,linecolor=blue]{->}(5.45,0)(5.1,0)
\psline[linewidth=1pt,linestyle=dotted,linecolor=blue]{-}(5.1,0)(4.9,0)
\psline[linewidth=.5pt,linecolor=blue]{->}(5.5,0)(6,0)
\psline[linewidth=.5pt,linecolor=blue]{->}(6.5,0)(7,0)
\psline[linewidth=1pt,linestyle=dotted,linecolor=blue]{-}(6.5,0)(6,0)
\psline[linewidth=.5pt,linecolor=blue]{->}(7.5,0)(7,0)
\psline[linewidth=1pt,linestyle=dotted,linecolor=blue]{-}(7.7,0)(9.5,0)
\psline[linewidth=.5pt,linecolor=blue]{->}(9.8,0)(10.5,0)
\psline[linewidth=.5pt,linecolor=blue]{<-}(11.5,0)(12,0)
\psline[linewidth=.5pt,linecolor=blue]{->}(10.9,0)(10.5,0)
\psline[linewidth=.5pt,linecolor=blue]{->}(11.45,0)(11.1,0)
\psline[linewidth=1pt,linestyle=dotted,linecolor=blue]{-}(11.1,0)(10.9,0)
\end{pspicture}
\end{center}
Thus the string $w(\ggz,\gz)$ is of the form
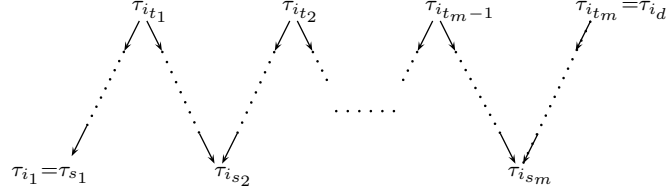
\begin{center}
\begin{pspicture}(-4.7,-.2)(3,2.2)
\put(-2,-.1){$_{\t_{i_{s_2}}}$}
\put(-1.1,2.1){$_{\t_{i_{t_2}}}$}
\put(-4.7,-.1){$_{\t_{i_1}=\t_{s_1}}$}
\put(1.9,-.1){$_{\t_{i_{s_m}}}$}
\put(.8,2.1){$_{\t_{i_{t_m-1}}}$}
\put(-3.1,2.1){$_{\t_{i_{t_1}}}$}

\psline[linewidth=1pt,linestyle=dotted]{-}(-.4,.7)(.4,.7)

\psline[linewidth=.5pt]{<-}(-1.2,1.5)(-1,1.9)
\psline[linewidth=.5pt]{->}(-1.7,.4)(-1.9,0)
\psline[linewidth=1pt,linestyle=dotted]{-}(-1.2,1.5)(-1.7,.4)

\psline[linewidth=1pt,linestyle=dotted]{-}(-.5,1.1)(-.7,1.5)
\psline[linewidth=.5pt]{->}(-.9,1.9)(-.7,1.5)

\psline[linewidth=.5pt]{->}(-2.9,1.9)(-2.7,1.5)
\psline[linewidth=.5pt]{<-}(-2,0)(-2.2,0.4)
\psline[linewidth=1pt,linestyle=dotted]{-}(-2.2,.4)(-2.7,1.5)

\psline[linewidth=.5pt]{->}(-3,1.9)(-3.2,1.5)%
\psline[linewidth=.5pt]{->}(-3.7,0.5)(-3.9,.1)
\psline[linewidth=1pt,linestyle=dotted]{-}(-3.7,0.5)(-3.2,1.5)

\psline[linewidth=1pt,linestyle=dotted]{-}(.5,1.1)(.7,1.5)
\psline[linewidth=.5pt]{->}(.9,1.9)(.7,1.5)

\psline[linewidth=.5pt]{<-}(1.2,1.5)(1,1.9)
\psline[linewidth=.5pt]{->}(1.8,.4)(2,0)
\psline[linewidth=1pt,linestyle=dotted]{-}(1.2,1.5)(1.8,.4)

\psline[linewidth=.5pt]{->}(2.3,.4)(2.1,0)
\psline[linewidth=1pt,linestyle=dotted]{-}(2.1,0)(3,1.9)
\psline[linewidth=.5pt]{<-}(2.8,1.5)(3,1.9)
\put(2.8,2.1){$_{\t_{i_{t_m}}=\t_{i_d}}$}
\end{pspicture}
\end{center}
\begin{itemize}
\item[$(1)$]
It is easy to see that top$M(\ggz,\gz)=M_I(\ggz,\gz)$ with $I=\{t_1, t_2, \ldots, t_m\}\subset\{1,2, \ldots, d\}$, and soc$M(\ggz,\gz)=M_J(\ggz,\gz)$ with $J=\{s_1, s_2,\ldots, s_{m}\}\subset\{1,2, \ldots, d\}$ where $1=s_1<t_1<s_2<t_2<\cdots <s_m<t_m=d$.

Therefore, the minimal injective resolution of $M(\ggz,\gz)$ is given by
$$0\lra M(\ggz,\gz)\lra \bigoplus\limits_{k=1}^m I(i_{s_k})\lra \bigoplus\limits_{k=1}^{m-1} I(i_{t_k})\oplus I(i_0)\oplus I(i_{d+1}),$$
where $I(i_k)$ is the indecomposable injective module corresponding to $\t_{i_k}$ and $I(i_k)=0$ if $\t_{i_k}$ is a boundary arc. By definition
$$\Ind_\ggz(\gz)=\deg(x_{i_0}x_{i_{d+1}}\prod_{k=1}^{m-1} x_{i_{t_k}}/\prod_{k=1}^m x_{i_{s_k}}).$$

On the other hand, we consider the unique $(\ggz,\gz)-$path $\az^0_\gz$ in Proposition \ref{uniqueness}:
$$\az_\gz^0=\az^0_1\az^0_2\cdots\az^0_{2d+1}=\t_{i_{0}}\underline{\t_{i_1}}\t_{i_2}\underline{\t_{i_2}}\t_{i_3}\cdots \t_{i_{t_1}}\underline{\t_{i_{t_1}}}\t_{i_{t_1}}\underline{\t_{i_{t_1+1}}}\t_{i_{t_1+1}}\cdots$$
$$\t_{i_{s_2-1}}\underline{\t_{i_{s_2}}}\t_{i_{s_2+1}}\underline{\t_{i_{s_2+1}}}\cdots\cdots \t_{i_{d-1}}\underline{\t_{i_{d-1}}}\t_{i_{d}}\underline{\t_{i_{d}}}\t_{i_{d+1}}$$
where the even positions are indicated by underlines. Then $y(\az^0_\gz)=1$ by definition and
$$x(\az^0_\gz)=\frac{\prod_{k\ odd}x_{\az^0_k}}{x_{i_1}x_{i_2}\cdots x_{i_d}}$$
$$=\frac{x_{i_0}x_{i_2}x_{i_3}\cdots x_{i_{t_1}}x_{i_{t_1}}x_{i_{t_1+1}}\cdots x_{i_{s_2-1}}x_{i_{s_2+1}}\cdots\cdots x_{i_{d-1}}x_{i_{d}}x_{i_{d+1}}}
{x_{i_1}x_{i_2}x_{i_3}\cdots x_{i_{t_1-1}}x_{i_{t_1}}x_{i_{t_1+1}}\cdots x_{i_{s_2-1}}x_{i_{s_2}}x_{i_{s_2+1}}\cdots\cdots x_{i_{d-1}}x_{i_{d}}}$$
$$=x_{i_0}x_{i_{d+1}}\prod_{k=1}^{m-1} x_{i_{t_k}}/\prod_{k=1}^m x_{i_{s_k}}.$$
This completes the proof.

\item[$(2)$] It suffices to prove that
$$\mbox{deg}(x(\az)y(\az))=\mbox{deg}(x(\az_\gz^0))$$
for each $\az\in C_\ggz(\gz).$ We know from section \ref{r-orient} that each complete $(\ggz,\gz)-$path $\az$ can be obtained from $\az^0_\gz$ by substituting pieces of segments (corresponding to substrings of $w(\ggz,\gz)$) in $\az_\gz^1$ for these in $\az_\gz^0$.

By the definition of $x(\az)$, it is sufficient to consider the complete $(\ggz,\gz)-$path obtained from $\az_\gz^0$ by substituting only one segment. Without loss of generality,
take $I=\{k,k+1,\ldots,t_1,t_1+1,\ldots,l\}\in S_\ggz(\gz)$, where $t_1< k\leq s_2\leq l<t_2.$ We consider the complete $(\ggz,\gz)-$path $$\az=\psi(I)=\az^0_1\underline{\az^0_2}\cdots \underline{\az^0_{2k-2}}(\t_{j_{k-1}}
\underline{\az^1_{2k}}\az^1_{2k+1}\cdots\az^1_{2t_1-1}\underline{\az^1_{2t_1}}\az^1_{2t_1+1}
$$
$$\cdots\underline{\az^1_{2l}} \t_{j_{l}})\underline{\az^0_{2l+1}}\cdots \az^0_{2t_2-1}\underline{\az^0_{2t_2}}\az^0_{2t_2+1}\cdots\cdots\underline{\az^0_{2d}} \az^0_{2d-1}$$
$$=\t_{i_{0}} \underline{\t_{i_1}} \t_{i_2} \underline{\t_{i_2}}\t_{i_3}\cdots
\t_{i_{t_1}}\underline{\t_{i_{t_1}}} \t_{i_{t_1}}\underline{\t_{i_{t_1+1}}}\t_{i_{t_1+1}}
\cdots\underline{\t_{i_{k-1}}}
(\t_{j_{k-1}}\underline{\t_{i_k}}\t_{i_{k+1}}\underline{\t_{i_{k+1}}}\t_{i_{k+2}}\cdots $$
$$\t_{i_{s_2}}\underline{\t_{i_{s_2}}}\t_{i_{s_2}}
\cdots\t_{i_{l-1}}\underline{\t_{i_l}}\t_{j_l})
\underline{\t_{i_{l+1}}}\t_{i_{l+2}}\cdots\t_{i_{s_2-1}}\underline{\t_{i_{s_2}}}\t_{i_{s_2+1}}\cdots\cdots
\t_{i_{d}}\underline{\t_{i_{d}}}\t_{i_{d+1}.}$$
Therefore
$$\mbox{deg}(x(\az)y(\az))-\mbox{deg}(x(\az_\gz^0))$$
$$=\mbox{deg}(x_{j_{k-1}}\cdot\prod_{a=k+1}^{s_2}x_{i_a}
\prod_{a=s_2}^{l-1}x_{i_a}\cdot x_{j_l}\cdot\prod_{a=k}^{l}y_{i_a}-\prod_{a=k-1}^{s_2-1}x_{i_a}
\prod_{a=s_2+1}^{l+1}x_{i_a})$$
$$=\mbox{deg}(\frac{x_{j_{k-1}}x_{i_{s_2}}x_{i_{s_2}}x_{j_l}}{{x_{i_{k-1}}x_{i_k}}x_{i_l}x_{i_{l+1}}})
+\mbox{deg}(\prod_{a=k}^ly_{i_a}),$$
where $x_{j_a}=1$ if $\t_{j_a}$ is a boundary arc.
By definition, we have
$$\mbox{deg}(\prod_{a=k}^ly_{i_a})=
\mbox{deg}(\prod_{a=k}^{s_2-1}\frac{x_{i_{a-1}}x_{j_{a}}}{x_{i_{a+1}}x_{j_{a-1}}}
\cdot
\frac{x_{i_{s_2-1}}x_{i_{s_2+1}}}{x_{j_{s_2-1}}x_{j_{s_2}}}\cdot
\prod_{a=s_2+1}^{l}
\frac{x_{i_{a+1}}x_{j_{a-1}}}{x_{i_{a-1}}x_{j_a}})$$
$$=\mbox{deg}(\frac{{x_{i_{k-1}}x_{i_k}}x_{i_l}x_{i_{l+1}}}{x_{j_{k-1}}x_{i_{s_2}}x_{i_{s_2}}x_{j_l}}).$$
Thus, $\mbox{deg}(x(\az)y(\az))=\mbox{deg}(x(\az_\gz^0))$ which completes the proof.
\end{itemize}
\end{pf}

From Proposition \ref{g-vector}(1), we can easily get the following corollary (Conjecture 8.3 in \cite{S}) which gives a recipe for computing the $\bg-$vector for any cluster variable $x^\ggz_\gz.$
\begin{cor} Let $\ggz$ be a triangulation of $(S,M)$, and $x^\ggz_\gz$ be the cluster variable associated to an internal arc $\gz$ in the cluster algebra ${\mathcal {A}}$, then
$$\bg^\ggz_\gz=\emph{Ind}_\ggz(\gz).$$
\end{cor}

\begin{rem} The above corollary was also proved in Proposition 6.2 in \cite{FK} for even more general cases, that is, 2-Calabi-Yau categories.
\end{rem}

\subsection{The expansion formula}

For each dimension vector $\be=(e_1,e_2,\ldots,e_n)^T,$ let

$$X^{\be}=\prod_{i=1}^n x_i^{e_i} \mbox{~and~} Y^{\be}=\prod_{i=1}^n y_i^{e_i}.$$
Hence,
$$y(\az)=\prod\limits_{\az_{2k}\ is\ \gz-oriented}y_{i_k}=Y^{\underline{\textbf{dim}}M_{I_\az}(\ggz,\gz)}$$
by definition for each complete $(\ggz,\gz)-$path $\az$. Define
$$\mu_\be(\ggz,\gz)=\sharp\{I\in S_\ggz(\gz)|\ \underline{\textbf{dim}}M_I(\ggz,\gz)=\be\}.$$ Then Schiffler's expansion formula can be interpreted in terms of representations as follows.

\begin{thm}\label{expansion}Let $\ggz$ be a triangulation of $(S,M)$, and $\gz$ be a curve in $(S,M)$, then
$$E^\ggz_\gz=\sum\limits_{\be}\mu_\be(\ggz,\gz)\ X^{\emph{Ind}_\ggz(\gz)+ B_\ggz\be}Y^{\be}.$$
\end{thm}
\begin{pf}Recall that for any two complete $(\ggz,\gz)-$paths $\az$ and $\bz$,
$$x(\az)y(\az)=x(\bz)y(\bz)\ \mbox{if\ and\ only\ if}\ y(\az)=y(\bz)$$ since
$E^\ggz_\gz$ is homogeneous by Proposition \ref{g-vector}.
Therefore $$E^\ggz_\gz=\sum\limits_{\az\in C(\gz)}x(\az)y(\az)$$
$$=\sum\limits_{I_\az\in S(\gz)}X^{\emph{Ind}_\ggz(\gz)+B_\ggz \underline{\textbf{dim}}M_{I_\az}(\ggz,\gz)}Y^{\underline{\textbf{dim}}M_{I_\az}(\ggz,\gz)}$$
$$=\sum\limits_{\be}\mu_\be(M(\ggz,\gz))\ X^{\emph{Ind}_\ggz(\gz)+B_\ggz \be}Y^{\be}.$$
As for the last equality, if $\az$ and $\bz$ are two complete $(\ggz,\gz)-$paths such that
$$\underline{\textbf{dim}}M_{I_\az}(\ggz,\gz)=\underline{\textbf{dim}}M_{I_\bz}(\ggz,\gz)=\be,$$
then both $I_\az$ and $I_\bz$ belong to $\mu_{\be}(\ggz,\gz)$.
\end{pf}

\begin{cor}Let $\ggz$ be a triangulation of $(S,M)$, and $x^\ggz_\gz$ be a cluster variable in the cluster algebra ${\mathcal{A}}_\ggz$ corresponding to an internal arc $\gz$ in $(S,M)$, then
 $$x^\ggz_\gz=\sum\limits_{\be}\mu_\be(M(\ggz,\gz))\ X^{\emph{Ind}_\ggz(\gz)+B_\ggz\be}Y^{\be}.$$
\end{cor}

\begin{example} Reconsider Example \ref{notation}, it is easy to see
$$B=B_\ggz=\left(
                             \begin{array}{ccc}
                               0 & -1 & 2 \\
                               1 & 0 & -1 \\
                               -2 & 1 & 0 \\
                             \end{array}
                           \right)$$
Hence $\deg(y_1)=(0,-1,2)^T$, $\deg(y_2)=(1,0,-1)^T$, $\deg(y_3)=(-2,1,0)^T.$
Take an internal arc $\gz$ as follows
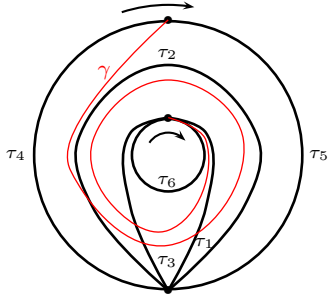
\begin{center}
\begin{pspicture}(-4,-1.8)(4,2)
\put(0,.5){\circle*{.1}}
\pscircle[linewidth=1pt](0,0){.5}
\pscircle[linewidth=1pt](0,0){1.8}
\put(0,1.8){\circle*{.1}}
\put(0,-1.8){\circle*{.1}}
\psarc[arcsepB=4pt]{<-}(0,0){.3}{-320}{-200}\psarc[arcsepB=4pt]{<-}(0,0){2}{80}{110}
\pscurve[linewidth=1pt]{-}(0,-1.8)(.5,-0.7)(.5,.3)(0,.5) \pscurve[linewidth=1pt]{-}(0,-1.8)(-.5,-0.7)(-.5,.3)(0,.5)
\pscurve[linewidth=1pt]{-}(0,-1.8)(1,-0.7)(1.2,.2)(0,1.2)(-1.2,.2)(-1,-0.7)(0,-1.8)
\pscurve[linewidth=.5pt,linecolor=red,showpoints=false]{-}(0,1.8)(-1.3,0.2)(-1.3,-0.2)(0,-1.2)(1,.2)(0,1)(-1,.2)(-.8,-.6)(0,-1)(.5,.2)(0,.5)
{\small\uput[u](0,1.1){$_{\t_2}$}\uput[l](-1.7,0){$_{\t_4}$}\uput[r](1.7,0){$_{\t_5}$}
\put(.35,-1.15){$_{\t_1}$}\put(-.15,-1.4){$_{\t_3}$}\uput[u](0,-.6){$_{\t_6}$}}
{\red\uput[l](-.7,1.1){$_\gz$}}
\end{pspicture}
\end{center}
then
$$w(\ggz,\gz)=\t_2\ra\t_3\ra\t_1\leftarrow\t_3.$$
$M(\ggz,\gz)=I_1$ is the indecomposable injective module over $J(Q_\ggz,W_\ggz)$ corresponding to the internal arc $\t_1$ and the minimal injective resolution of $M(\ggz,\gz)$ is given by
$$ 0\lra M(\ggz,\gz)\lra I_1\lra 0,$$
hence $\Ind_\ggz(\gz)=(-1,0,0)^T$ by definition. Therefore Theorem \ref{expansion} implies that
$$E^\ggz_\gz=\mu_{(0,0,0)}(\ggz,\gz)\cdot X^{(-1,0,0)}\cdot Y^{(0,0,0)}$$$$
+\mu_{(1,0,0)}(\ggz,\gz)\cdot X^{(-1,1,-2)}\cdot Y^{(1,0,0)}$$
$$+\mu_{(1,0,1)}(\ggz,\gz)\cdot X^{(1,0,-2)}\cdot Y^{(1,0,1)}$$$$
+\mu_{(1,1,1)}(\ggz,\gz)\cdot X^{(0,0,-1)}\cdot Y^{(1,1,1)}$$
$$+\mu_{(1,0,2)}(\ggz,\gz)\cdot X^{(3,-1,-2)}\cdot Y^{(1,0,2)}$$$$
+\mu_{(1,1,2)}(\ggz,\gz)\cdot X^{(2,-1,-1)}\cdot Y^{(1,1,2)}$$
$$=\frac{1}{x_1x_2x_3^2}(x_2x_3^2+x_2^2y_1+2x_1^2x_2y_1y_3$$$$+x_1x_2x_3y_1y_2y_3
+x_1^4y_1y_3^2+x_1^3x_3y_1y_2y_3^2).
$$
\end{example}

\section{Caldero-Chapoton map and Schiffler's expansion formula}
We show in this section that the Caldero-Chapoton map (or cluster character) and Schiffler's expansion formula coincide.
\subsection{Caldero-Chapoton map}

Let $\mc$ be a Hom-finite, Krull-Schmidt, 2-Calabi-Yau triangulated $k-$category (where $k$ is an algebraically closed field), and assume that $\mc$ admits a cluster-tilting object $T=T_1\oplus\cdots\oplus T_n$ with $n$ direct summands. Let $A=\End^{op}_\mc(T)$ be the endomorphism algebra of $T$ in $\mc$. Recall that there is an equivalence of categories:
$$\Ext_\mc^1(T,?): \mc/T\lra {\bf mod}A.$$

Given a finite dimensional $A-$module $M$, and a dimension vector $\be\in {\mathbb{N}}^n$, we denote by
$$Gr_{\be}(M)=\{N\subseteq M\ \mbox{submodule~such~that}\ \underline{dim}N=\be\}$$
the variety of $\be-$dimensional submodules of $M$. Note that $Gr_{\be}(M)$ is a projective variety (a subvariety of the classical Grassmannian of subspaces of $M$). We denote by $\chi(Gr_\be(M))$ its Euler-Poincar\'e characteristic.

\begin{rem} \label{grothedick} Let ${\mathcal {T}}=\mbox{add}T$ be the full subcategory of $\mc$ whose objects are all direct factors of direct sums of $T,$ then ${\mathcal {T}}$ is equivalent to the category of finitely generated projective (or injective) modules over $A.$
\end{rem}
We denote by ${\mathcal {K}}_0({\mathcal {T}})$ the Grothendieck group of the additive category ${\mathcal {T}}.$

\begin{lem}[\cite{KR}]\label{approximation} For each object $M$ of $\mc$, there exists a triangle
$$T_1\lra T_0\lra M\lra T_1[1]$$
where $T_0$, $T_1$ belong to ${\mathcal {T}}$. Moreover, as an element in ${\mathcal {K}}_0({\mathcal {T}})$, $$[T_0]-[T_1]$$ does not depend on the choice of this triangle.
\end{lem}
We define the {\bf index} $\mbox{Ind}_T(M)$ with respect to $T$ of an object $M$ in $\mc$ as the element $[T_0]-[T_1]$ of
${\mathcal {K}}_0({\mathcal {T}})$.
\medskip

The {\em Caldero-Chapoton map} (or {\em Cluster character} in sense of \cite{Palu}) with respect to $T$
$$X^T_?: \Ind \mc\lra {\mathbb{Q}}(x_1,\ldots,x_n)$$
is defined as
$$X^T_M:=\sum_{\be}\chi(Gr_{\be}(\Ext^1_\mc(T,M)))\cdot{X}^{\Ind_T(M)+B_T\be}$$
where $B_T$ is the skew-symmetric matrix induced by the quiver of $\End_\mc(T)$.

Let ${\mathcal {A}}^0_T$ be the cluster algebra with initial seed $(B_T, {\bf x},1)$ where  ${\bf x}=(x_1,\ldots,x_n).$

\begin{prop}[\cite{Palu}]\label{palu} The cluster character $X^T_?$ on $\mc$ induces a bijection between the set of reachable rigid indecomposable objects of $\mc$ to the set of cluster variables in ${\mathcal {A}}^0_T.$
\end{prop}

\subsection{Cluster category of a marked surface}

Let $(S,M)$ be a marked surface without punctures, recall that the cluster category $\mc_{(S,M)}$ of the
surface is defined in \cite{A,BZ}. In fact, a quiver with potential $(Q_\ggz,W_\ggz)$ whose Jacobian algebra is finite-dimensional is associated to each triangulation $\ggz$ of $(S,M)$, and the cluster category $\mc_{(S,M)}$ is just defined as the cluster category of the quiver with potential $(Q_\ggz^{op},W_\ggz^{op})$ introduced by Amiot (we refer \cite{A} for more details).
Note that $\mc_{(S,M)}$ is a Hom-finite, Krull-Schmidt, 2-Calabi-Yau triangulated category and it does not depend on the triangulation of $(S,M)$ up to triangle equivalence.

It has been show in \cite{BZ} that each triangulation $\ggz$ of $(S,M)$ yields a cluster-tilting object $T_\ggz$ in $\mc_{(S,M)}$ and  this yields an equivalence of categories by \cite{KZ}
$$\Ext^1_{\mc_{(S,M)}}(T_\ggz,?):~~ \mc_{(S,M)}/T_\ggz \overset{\sim}{\lra} {\bf mod}J(Q_\ggz,W_\ggz).$$
Based on the above equivalence, a geometric characterization of the indecomposable objects in $\mc_{(S,M)}$
has been given in \cite{BZ}. Note that in the special case where $S$ is a disc, the characterization was given in \cite{CCS}.

\begin{thm}[\cite{CCS,BZ}]\label{indec}A parametrization of the isoclasses of indecomposable objects in $\mc_{(S,M)}$ is given by ``string objects" and ``band objects", where
\begin{itemize}
\item[$(1)$] the string objects are indexed by the homotopy classes of non-contractible curves in $(S,M)$ which are not homotopic to a boundary segment of $(S,M)$.
\item[$(2)$] the band objects are indexed by $k^* \times {\Pi^*_1 (S,M)}/ \sim$, where $k^*=k\setminus\{0\}$ and $\Pi^*_1 (S,M)/\sim$  is given by the non-trivial elements of the fundamental group of $(S,M)$ subject to the equivalence relation generated by $a\sim a^{-1}$ and cyclic permutation.
\end{itemize}
\end{thm}
\begin{rem}Let $\gz$ be a curve in $(S,M)$, if we view $\gz$ as a string object in $\mc_{(S,M}$, then
$$\Ext^1_{\mc_{(S,M)}}(T_\ggz,\gz)=M(\ggz,\gz).$$
\end{rem}

Theorem \ref{indec} and Proposition \ref{palu} imply that
the cluster character $X^\ggz_?:=X^{T_\ggz}_?$ induces a map from the curves in $(S,M)$ to ${\mathbb{Q}}(x_1,\ldots,x_n)$. On the other hand, Schiffler-Thomas' formula $E^\ggz_?(x_1,\ldots,x_n,1,\ldots,1)$ (see Remark \ref{Schiffler-Thomas' formula}) also produces a map from the curves in $(S,M)$ to ${\mathbb{Q}}(x_1,\ldots,x_n)$. These two maps are expected to be the same.

\subsection[The Caldero-Chapoton map and Schiffler-Thomas' formula]{Caldero-Chapoton map and Schiffler-Thomas' formula}
We show in the this subsection that the cluster character $X^\ggz_?$ of $\mc_{(S,M)}$ of with respect $T_\ggz$ and Schiffler-Thomas' expansion formula $E^\ggz_?$ coincide for curves in $(S,M)$, that is, the ``string objects'' in $\mc_{(S,M)}.$ The following lemma is important in this subsection.

\begin{lem}\label{index-lemma}Let $\gz$ be a curve in $(S,M)$. If there is a short exact sequence in {\bf mod}$J(Q_\ggz,W_\ggz)$:
$$0\lra M(\ggz,\gz)\overset{f}{\lra} Q_0\overset{g}{\lra} Q_1\lra N$$
with $Q_0,Q_1$ injective and $N$ indecomposable, then
$${\emph{Ind}}_\ggz(\gz)=[Q_1]-[Q_0]\in {\mathcal {K}}_0(\emph{add}{{\mathcal {I}}})$$
where ${\mathcal {I}}$ is the category of finitely generated  injective modules over $J(Q_\ggz,W_\ggz).$
\end{lem}
\begin{pf}We take the minimal injective resolution of $M(\ggz,\gz):$
  $$0\lra M(\ggz,\gz)\overset{f'}{\lra} I_0\overset{g'}{\lra} I_1.$$
By definition, $(I_0, f')$ is the injective envelope of $M(\ggz,\gz)$, hence
$$Q_0\simeq I_0\oplus I$$
with $I$ injective. Consider the following commutative diagram:
 \begin{center}
\begin{pspicture}(-3.4,-.7)(5.7,1.2)
\put(-3.8,.7){$0$}
\put(-3.8,-.7){$0$}
\psline[linewidth=.7pt]{->}(-3.4,0.8)(-2.1,.8)
\psline[linewidth=.7pt]{->}(-3.4,-0.6)(-2.1,-.6)
\put(-2,.7){$M(\ggz,\gz)$}
\put(-2,-.7){$M(\ggz,\gz)$}
\psline[linewidth=.7pt]{->}(-.5,0.8)(.8,.8)
\psline[linewidth=.7pt]{->}(-.5,-0.6)(.8,-.6)
\psline[linewidth=.7pt]{<-}(-1.7,-0.3)(-1.7,.5)
\psline[linewidth=.7pt]{<-}(1.2,-0.3)(1.2,.5)

\put(1,.7){$I_0$}
\put(1,-.7){$Q_0$}
\psline[linewidth=.7pt]{->}(1.5,0.8)(2.8,.8)
\psline[linewidth=.7pt]{->}(1.5,-0.6)(2.8,-.6)
\put(3,.7){$\mbox{coker}f'$}
\put(3,-.7){$\mbox{coker}f$}
\psline[linewidth=.7pt]{->}(4.3,0.8)(5.6,.8)
\psline[linewidth=.7pt]{->}(4.3,-0.6)(5.6,-.6)
\put(5.7,.7){$0$}
\put(5.7,-.7){$0$}
\psline[linewidth=.7pt]{<-}(3.4,-0.3)(3.4,.5)

\uput[u](0,.7){$_{f'}$}
\uput[u](0,-.7){$_{f}$}
\uput[l](-1.7,0.1){$_{1}$}
\end{pspicture}
\end{center}

 Schanuel's Lemma implies that $$\mbox{coker}f'\oplus Q_0\simeq \mbox{coker}f'\oplus I_0\oplus I\simeq  \mbox{coker}f\oplus I_0.$$
Therefore $\mbox{coker}f'\oplus I \simeq \mbox{coker}f.$  By the definition of minimal injective resolution,
$(I_1, i')$ is the injective envelope of $\mbox{coker}f'\simeq \mbox{im}g',$ see the following commutative diagram
\begin{center}
\begin{pspicture}(-3.4,-.7)(5.7,1.2)
\put(-3.8,.7){$0$}
\put(-3.8,-.7){$0$}
\psline[linewidth=.7pt]{->}(-3.4,0.8)(-2.1,.8)
\psline[linewidth=.7pt]{->}(-3.4,-0.6)(-2.1,-.6)
\put(-2,.7){$\mbox{coker}f'$}
\put(-2,-.7){$\mbox{coker}f$}
\psline[linewidth=.7pt]{->}(-.5,0.8)(.8,.8)
\psline[linewidth=.7pt]{->}(-.5,-0.6)(.8,-.6)
\psline[linewidth=.7pt]{<-}(-1.7,-0.3)(-1.7,.5)
\psline[linewidth=.7pt]{<-}(1.2,-0.3)(1.2,.5)

\put(1,.7){$I_1$}
\put(1,-.7){$Q_1$}
\psline[linewidth=.7pt]{->}(1.5,0.8)(2.8,.8)
\psline[linewidth=.7pt]{->}(1.5,-0.6)(2.8,-.6)
\put(3,.7){$\mbox{coker}i'$}
\put(3,-.7){$\mbox{coker}i$}
\psline[linewidth=.7pt]{->}(4.3,0.8)(5.6,.8)
\psline[linewidth=.7pt]{->}(4.3,-0.6)(5.6,-.6)
\put(5.7,.7){$0$}
\put(5.7,-.7){$0$}
\psline[linewidth=.7pt]{<-}(3.4,-0.3)(3.4,.5)

\uput[u](0,.7){$_{i'}$}
\uput[u](0,-.7){$_{i}$}
\uput[l](-1.7,0.1){$_{h}$}
\end{pspicture}
\end{center}
where $h: \mbox{coker}f'\lra\mbox{coker}f\simeq\mbox{coker}f'\oplus I$ is an embedding.
We claim that
$(Q_1,i)$ is also the injective envelope of $\mbox{coker}f\simeq \mbox{coker}f'\oplus I .$ Otherwise,
we can write $Q_1\simeq I_1\oplus I\oplus I'$ with $I'\neq 0$ an injective module. It implies that
$$\mbox{coker}g\simeq\mbox{coker}i\simeq  Q_1/\mbox{im}i\simeq  (I_1\oplus I\oplus I')/(\mbox{im}i'\oplus I)
\simeq \mbox{coker}i'\oplus I'. $$
This contradicts the fact that $N$ is indecomposable since $\mbox{coker}g$ which is isomorphic to a submodule of $N$ contains a non-zero injective summand.
Therefore
$$Q_1\simeq I_1\oplus I$$
which implies that
$$\mbox{Ind}_\ggz(\gz)=[I_1]-[I_0]=[Q_1]-[Q_0].$$
\end{pf}

\begin{thm}Let $\gz$ be a curve in $(S,M)$, then
$$X^\ggz_\gz=E^\ggz_\gz(x_1,\ldots,x_n,1,\ldots,1).$$
\end{thm}
\begin{pf}It suffices to consider a curve $\gz\not\in\ggz$. By Lemma \ref{approximation} there exists a triangle
$$\t_{a_1}\overset{f}{\lra} \t_{a_0}\overset{g}{\lra} \gz\overset{h}{\lra}  \t_{a_1}[1]$$
where ${\mathcal {T}}=\mbox{add}T_\ggz$, $\t_{a_1}$ and $\t_{a_2}$ belong to add${\mathcal {T}}$.
By applying $\Ext^1_\mc(T_\ggz,?)$ to the above triangle, there is a short exact sequence of $J(Q_\ggz,W_\ggz)-$modules
$$0\lra M(\ggz,\gz)\overset{Ext_\mc^1(T_\ggz,h)}{\lra} \Hom_\mc(T_\ggz,\t_{a_1}[2])\overset{Ext_\mc^1(T_\ggz,f[1])}{\lra} $$
$$\Hom_\mc(T_\ggz,\t_{a_0}[2])\overset{Ext_\mc^1(T_\ggz,g[1])}{\lra} \Hom_\mc^1(T_\ggz,\gz[2]).$$
Since $\mc$ is a 2-Calabi-Yau $k-$category, the exact sequence becomes
$$0\lra M(\ggz,\gz)\overset{Ext_\mc^1(T_\ggz,h)}{\lra} D\Hom_\mc(\t_{a_1},T_\ggz)\overset{Ext_\mc^1(T_\ggz,f[1])}{\longrightarrow} $$
$$D\Hom_\mc(\t_{a_0},T_\ggz)\overset{Ext_\mc^1(T_\ggz,g[1])}{\lra} \tau(M(\ggz,\gz))$$
where $D\Hom_\mc(\t_{a_1},T_\ggz)$ and $D\Hom_\mc(\t_{a_0},T_\ggz)$ are injective modules over $J(Q_\ggz,W_\ggz)$. Since $\tau(M(\ggz,\gz))$ is indecomposable or 0 in {\bf mod}$J(Q_\ggz,W_\ggz)$, Lemma \ref{index-lemma} and Remark \ref{grothedick} implies
$$\Ind_\ggz(\gz)=\Ind_{T_\ggz}(\gz)=[\t_{a_0}]-[\t_{a_1}]\in {\mathcal {K}}_0({\mathcal {T}}).$$

Moreover, Poettering in \cite{Poettering} gives a recipe to compute the Euler characteristic of a string modules which implies that
$$\mu_{\be}(\ggz,\gz)=\chi(Gr_{\be}(M(\ggz,\gz))).$$

Therefore Theorem \ref{expansion} shows
$$E^\ggz_\gz(x_1,\ldots,x_n,1,\ldots,1)=\sum\limits_{\be}\mu_\be(\ggz,\gz)\ X^{\emph{Ind}_\ggz(\gz)+B_\ggz\be}$$
$$=\sum_{\be}\chi(Gr_{\be}(\Ext^1_\mc(T_\ggz,\gz)))\cdot{X}^{\emph{Ind}_{T_\ggz}(\gz)+B_{T_\ggz}\be}=
X^\ggz_\gz.$$
Note that $B_\ggz=B_{T_\ggz}$ which completes the proof.
\end{pf}

\end{document}